# A superprocess model for eDNA transport in random river environment

Hidekazu Yoshioka[1, *], Motoh Tsujimura[2], Ayumi Hashiguchi[3]

[1] Graduate School of Advanced Science and Technology, Japan Advanced Institute of Science and Technology, Japan

[2] Graduate School of Commerce, Doshisha University, Japan

[3] Faculty of Environmental, Life, Natural Science and Technology, Okayama University, Japan

* E-mail: yoshih@jaist.ac.jp

ORCID: https://orcid.org/0000-0002-5293-3246 (Yoshioka), https://orcid.org/0000-0001-6975-9304 (Tsujimura)

**Abstract**

Environmental DNA (eDNA) is increasingly used to monitor aquatic organisms in rivers, yet predicting its concentration and assessing the reliability of such predictions remain challenging because river flow fluctuates with long memory and models of eDNA dynamics are inevitably misspecified. We propose a stochastic modeling framework for eDNA concentration in rivers that addresses both issues. The eDNA concentration is described by a stochastic partial differential equation with transport, decay, source, and square-root multiplicative noise, and the flow velocity is modeled by a superposition of Ornstein–Uhlenbeck-type processes with long memory. We show that the long memory of the flow velocity induces superdiffusive transport of eDNA, with the variance of the transported distance growing superlinearly in the ensemble sense. The proposed model admits a closed-form Laplace functional, based on which we prove the existence and uniqueness in law of the solution on an unbounded domain. We further evaluate the effects of model uncertainty through a robust control formulation based on the relative entropy, in which the lower and upper bounds of the expected eDNA concentration, together with the worst-case distortions, are obtained in closed form. The framework is applied to eDNA and flow velocity data collected for the upstream migration of *Plecoglossus altivelis altivelis* in a river in Japan. The results demonstrate that the closed-form bounds provide a computationally inexpensive tool for quantifying the range of eDNA concentrations that can be expected under model uncertainty.

***Acknowledgments*** We would like to express our sincere gratitude to the members of the Hii River Fisheries Cooperative for their tremendous cooperation during the field survey. The eDNA analysis was outsourced to the Environmental Pollution Research Center Co., Ltd., Japan.

## 1. Introduction

Environmental DNA (eDNA), the DNA released by organisms into water, has become a widely used tool for monitoring aquatic species because it requires only water sampling and poses minimal disturbance to ecosystems (Beng and Corlett 2020; Rees et al. 2014). However, eDNA concentrations fluctuate significantly in space and time (Augustine et al. 2026a; Inui et al. 2021; Yuan et al. 2026), owing not only to the population dynamics of the target species but also to the uncertainty in the elementary processes of shedding, decay, and transport (Lamb et al. 2022). Rivers deserve particular attention in this respect: they serve as habitats and migration routes of many aquatic species, and their flow, driven by hydrological events such as rainfall, fluctuates stochastically and often with long memory (Calvani and Perona 2025; Houenafa et al. 2025; Studnicka et al. 2026). Predicting eDNA concentrations in rivers and quantifying the reliability of such predictions therefore require models that account for both the randomness of river flow and the possible misspecification of the model itself.

Existing models of riverine eDNA transport are mostly deterministic, based on linear partial differential equations (Augustine et al. 2026b; Nukazawa et al. 2018; Xiong et al. 2025) or mass balance on river networks (Wu et al. 2025). As a stochastic alternative, Yoshioka (2026a) proposed a stochastic partial differential equation (SPDE) for the eDNA concentration of upstream-migrating fish in a bounded river reach, but the flow velocity and decay rate were assumed to be constant, and model uncertainty was not addressed. Meanwhile, passive particles transported by a long-memory velocity field are known to exhibit superdiffusion (Grahovac and Kevei 2025), suggesting that the randomness of river flow can substantially affect eDNA transport; representing this effect by a classical or fractional diffusion term is either inconsistent with the long memory or unrealistic for rivers with a dominant flow direction (Monroy and Raposo 2024).

The aim of this paper is to develop a stochastic framework for eDNA concentration in rivers that accounts for long-memory flow and model uncertainty while remaining analytically tractable. Our contributions are threefold. First, we formulate the eDNA concentration as a superprocess, namely an SPDE with transport, decay, source, and square-root multiplicative noise (Etheridge 2000), driven by a flow velocity modeled as a superposition of Ornstein–Uhlenbeck processes with long memory (Barndorff-Nielsen 2001; Yoshioka 2026c), and we show that this randomization induces superdiffusive transport of eDNA in the ensemble sense. Second, we prove that the model is well posed on an unbounded domain and admits a closed-form Laplace functional. Third, following the robust control framework of Hansen and Sargent (2001), we quantify the effects of model uncertainty through the relative entropy and derive the lower and upper bounds of the expected eDNA concentration, together with the worst-case distortions, in closed form; these bounds can be interpreted as certainty-equivalent risk measures. To our knowledge, model uncertainty of this type has been studied for SPDEs only with additive noise (Baltas et al. 2019) and not for superprocesses. We apply the framework to eDNA and flow velocity data collected for the upstream migration of *Plecoglossus altivelis altivelis* in a river in Japan.

The rest of this paper is organized as follows. **Section 2** presents the model and its mathematical analysis. **Section 3** addresses model uncertainty. **Section 4** presents the application. Section 5 concludes the paper. Appendices contain proofs, auxiliary results and eDNA data.

## 2. Mathematical model

Throughout this paper, we use a complete filtered probability space $\left(\Omega,\mathcal{F},\left(\mathcal{F}_t\right)_{t\geq 0},\mathbb{P}\right)$ as a standard technique in stochastic analysis, and we denote expectation by $\mathbb{E}$. Often, subscripts and arguments of variables and coefficients are omitted when there is no confusion. The collection of all functions $\phi:\mathbb{R}\to\mathbb{R}$ that are continuously differentiable arbitrary times is denoted by $C^{\infty}$. The collection of all functions $\phi:\mathbb{R}\to\mathbb{R}$ that are bounded and continuous is denoted by $C_b$. We set $C_b^{\infty}=C^{\infty}\bigcap C_b$. The collection of all nonnegative finite measures on $\mathbb{R}$ is denoted by $\mathcal{M}_+$. For any $T>0$, the collection of all time-dependent nonnegative measures that are càdlàg in time endowed with the Skorohod topology is denoted by $\mathbb{D}\left(\left[0,T\right],\mathcal{M}_+\right)$.

### 2.1 Superprocess model of eDNA concentration

We consider eDNA dynamics along a river as a one-dimensional space. Unlike oceans or lakes, rivers can be conceptualized as one-dimensional domains, along which the migration of aquatic species and eDNA transport are constrained.

The eDNA concentration at time $t$ and location $x$ is governed by the following SPDE:

$$\underbrace{\frac{\partial X_t\left(x\right)}{\partial t}}_{\text{Temporal change}}=\underbrace{u_t\frac{\partial X_t\left(x\right)}{\partial x}}_{\text{Advection}}-\underbrace{R_t X_t\left(x\right)}_{\text{Decay}}+\underbrace{f_t\left(x\right)}_{\text{Source}}+\underbrace{\sigma_t\sqrt{X_t\left(x\right)}\dot{W}}_{\text{Noise}}, \tag{1}$$

subject to an initial condition $X_0\left(x\right)$ ($x\in\mathbb{R}$), where $X_t\left(x\right)$ is the eDNA concentration of an aquatic species, $u_t>0$ is the flow velocity (the flow direction is toward negative $x$), $R_t>0$ is the decay rate, $f_t\left(x\right)>0$ is the source and hence shedding, $\sigma_t>0$ is the noise intensity, and $\dot{W}$ is a space-time white noise. The essential difference between the SPDE of Yoshioka (2026a) and the present model (1) is that the flow velocity was constant in the former, whereas its environmental, stochastic dependence is the focus in the latter; instead, the source in the proposed model is made simpler. The use of white noise is consistent with the finding that temporal residuals in eDNA concentration data have short-range correlations (Augustine et al. 2026b), which can be negligible for problems with longer timescales, such as monthly ones.

The eDNA concentration $X_t\left(x\right)$ is considered a superprocess because it depends both on space and time and is possibly measure-valued; the SPDE (1) contains the non-Lipschitz multiplicative noise term $\sigma_t\sqrt{X_t\left(x\right)}\dot{W}$, which ensures the nonnegativity of the solution in the sense of measures, a common feature of superprocess models (Etheridge 2000). Moreover, this noise term leads to a closed-form Laplace functional. For example, if there is no source and the parameters are constants, then we obtain the closed-form formula

$$\mathbb{E}\left[\exp\left(-\int_{\mathbb{R}}\phi(x)X_t(x)\mathrm{d}x\right)\right]=\exp\left(-\int_{\mathbb{R}}\frac{e^{-Rt}\phi(x-ut)}{1+\frac{\sigma^2}{2R}\left(1-e^{-Rt}\right)\phi(x-ut)}X_0(x)\mathrm{d}x\right) \tag{2}$$

for any sufficiently regular functions $\phi:\mathbb{R}\to\mathbb{R}$. This remarkable property carries over to the more complex cases studied in this paper.

We consider eDNA dynamics during a time interval $(0,T)$ with some fixed $T>0$. The SPDE (1) is formal because the "solution" $X$ may be a measure. We therefore consider the weak formulation:

$$\begin{aligned}&\underbrace{\int_{\mathbb{R}}\phi(x)X_t(\mathrm{d}x)-\int_{\mathbb{R}}\phi(x)X_0(\mathrm{d}x)}_{\text{Increment}}\\&=-\underbrace{\int_0^t\int_{\mathbb{R}}u_s\frac{\mathrm{d}\phi(x)}{\mathrm{d}x}X_s(\mathrm{d}x)\mathrm{d}s}_{\text{Advection}}-\underbrace{\int_0^t\int_{\mathbb{R}}R_s\phi(x)X_s(\mathrm{d}x)\mathrm{d}s}_{\text{Decay}}+\underbrace{\int_0^t\int_{\mathbb{R}}\phi(x)f_s(x)\mathrm{d}x\mathrm{d}s}_{\text{Source}}+\underbrace{\hat{M}_t(\phi)}_{\text{Noise}}\end{aligned} \tag{3}$$

for any $\phi\in C_b^\infty$, where $\hat{M}_t(\phi)$ is a martingale with quadratic variation

$$\left\langle\hat{M}_t(\phi)\right\rangle=\int_0^t\int_{\mathbb{R}}\sigma_s^2(\phi(x))^2X_s(\mathrm{d}x)\mathrm{d}s\,, \tag{4}$$

provided that $X\in\mathbb{D}\left([0,T],\mathcal{M}_+\right)$. We assume $X_0\in\mathcal{M}_+$. The weak formulation (3) is obtained by multiplying $\phi$ by (1) and then taking integration by parts.

The advantages of the weak formulation (3) are threefold. First, we can address measure-valued processes for which the spatial partial derivative $\frac{\partial X}{\partial x}$ does not necessarily exist in the classical sense. Second, we do not have to directly deal with the noise term driven by white noise; indeed, if $X$ is nonnegative and considered a function, the last term in (3) may be understood as follows:

$$\hat{M}_t(\phi)=\int_0^t\int_{\mathbb{R}}\sigma_s\phi(x)\sqrt{X_s(x)}W(\mathrm{d}s,\mathrm{d}x)\,, \tag{5}$$

where $W$ is space-time white noise (Definition 1.3.1 in Dalang and Sanz-Solé (2026)). The interpretation of the square root $\sqrt{X_s(x)}$ becomes difficult when $X$ is a measure. The weak formulation (3) avoids this theoretical issue because it only indirectly specifies $\hat{M}_t(\phi)$ through its quadratic variation. Third, the weak formulation can rigorously handle a wide range of eDNA profiles that are possibly very irregular (neither continuous nor differentiable). Phenomenologically, this noise term can be interpreted as stemming from unresolved factors such as stochasticity in each elementary processes, uneven mixing due to turbulence, and the process by which eDNA particles break down into smaller fragments with this fragmentation viewed as branching.

***Remark 1.*** In application, the spatial domain can be replaced by a bounded interval. In this case, eDNA will be freely discharged from downstream and a free-outflow condition should be specified at the downstream-end, whereas no boundary condition is necessary at the upstream-end.

***Remark 2.*** The most famous superprocesses are super-Brownian motion and its variants, which are essentially stochastic heat equations with non-Lipschitz multiplicative noise terms (Dawson et al. 1989; Le Doussal 2022). Population dynamics models as superprocesses with age structures have also been investigated by focusing on occupation measures (Cheng 2024) and perturbation analysis (Méléard and Tran 2012). Superprocess models have also been investigated in view of infinite-dimensional affine and polynomial SDEs (Schmidt et al. 2020), whose Laplace functionals or moment equations are obtained analytically, and an SPDE related to (1) has been studied computationally in the context of energy finance by assuming the existence of a solution as a function (Cuchiero et al. 2025). Spatially one-dimensional ecological models based on superprocesses (Bhat 2025; Week et al. 2021) and SPDEs with advection and white noise for soil water storage (Chang et al. 2025) have also been studied, but unlike this study, they did not address applications to eDNA dynamics.

## 2.2 supOU process of flow velocity

The flow velocity data obtained at our study site suggest that $u_t$ is a positive long-memory process driven by jumps (**Section 4.3**). We assume that the flow velocity $u$ is a stationary supOU process, which is a long-memory process obtained by suitably integrating independent Ornstein–Uhlenbeck processes with different timescales (Barndorff-Nielsen 2001; Yoshioka 2026c):

$$u_t = \int_{r=0}^{r=+\infty} \int_{z=0}^{z=+\infty} \int_{s=-\infty}^{s=t} z e^{-r(t-s)} \mathfrak{N}(\mathrm{d}r, \mathrm{d}z, \mathrm{d}s), \quad t \in \mathbb{R}, \tag{6}$$

where $r > 0$ is the reversion speed, $z > 0$ is the distributed jump size, $\mathfrak{N}$ is a Poisson random measure on $(0,+\infty) \times (0,+\infty) \times \mathbb{R}$ whose compensation measure equals $\pi(\mathrm{d}r)\nu(\mathrm{d}z)\mathrm{d}s$, $\pi(\mathrm{d}r)$ is a probability measure such that $\int_0^{+\infty} r^{-1}\pi(\mathrm{d}r) < +\infty$ (necessary to have finite moments of $u_t$, e.g., Yoshioka (2026c)), called the reversion measure, of $r$, and $\nu(\mathrm{d}z)$ is a Lévy measure with $\int_0^{+\infty} z^m \nu(\mathrm{d}z) \in (0,+\infty)$ ($m \in \mathbb{N}$) and does not have singularity at 0; this regularity assumption on $\nu(\mathrm{d}z)$ ensures the existence of all positive-order moments of flow velocity $u_t$. We assume that the Poisson random measure $\mathfrak{N}$ is independent of the driving white noise of the superprocess model.

We assume that the reversion measure $\pi$ has a gamma density $p(r) = r^{\alpha-1}e^{-r/\beta} / \left(\beta^\alpha \Gamma(\alpha)\right)$, where $\alpha > 1$ is the shape parameter, $\beta > 0$ is the scale parameter, and $\Gamma$ is the gamma function. The autocorrelation function $\rho$ of $u$ is given by (Yoshioka 2026c)

$$\rho(\tau) = \frac{1}{(1+\beta\tau)^{\alpha-1}}, \tag{7}$$

where $\tau \geq 0$ is time lag. The process $u$ admits long memory if $\alpha \in (1,2]$, yielding $\int_0^{+\infty} \rho(\tau)\,\mathrm{d}\tau = +\infty$. Moments of $u$ are obtained analytically as well; as the name implies, the supOU process is formed by superimposing OU processes that are linear; see Yoshioka (2026c).

Theorem 3.12 in Barndorff-Nielsen and Stelzer (2011) shows that our supOU process admits paths that a.s. have bounded variations in each compact time interval because $\int_0^{+\infty} z\nu(\mathrm{d}z) \in (0,+\infty)$. Additionally, this process is bounded in each compact time interval because of the definition of the bounded variation functions, e.g., Theorem 17.3(i) in Laczkovich and Sós (2017). The supOU process is a.s. positive for each compact time interval (Barndorff-Nielsen 2001). Indeed, for each $t \in [0,T]$, we have

$$\begin{aligned} u_t &= \int_{r=0}^{r=+\infty} \int_{z=0}^{z=+\infty} \int_{s=-\infty}^{s=0} z e^{-r(t-s)} \mathfrak{N}(\mathrm{d}r,\mathrm{d}z,\mathrm{d}s) + \int_{r=0}^{r=+\infty} \int_{z=0}^{z=+\infty} \int_{s=0}^{s=t} z e^{-r(t-s)} \mathfrak{N}(\mathrm{d}r,\mathrm{d}z,\mathrm{d}s) \\ &\geq \int_{r=0}^{r=+\infty} \int_{z=0}^{z=+\infty} \int_{s=-\infty}^{s=0} z e^{-r(T-s)} \mathfrak{N}(\mathrm{d}r,\mathrm{d}z,\mathrm{d}s), \end{aligned} \tag{8}$$

where the last line is positive, implying that each sample path of the supOU process on $[0,T]$ has a (pathwise) minimum value that is positive.

These properties of $u$ imply that the following integrated process, which is the flow distance of a passive particle during $(0,t)$, is well defined:

$$l_t = \int_0^t u_s \mathrm{d}s, \quad t \geq 0. \tag{9}$$

Because $u$ is stationary, the average $\mathbb{E}[l_t]$ of $l_t$ is obtained as

$$\mathbb{E}[l_t] = \mathbb{E}[u_0]t, \quad t \geq 0 \tag{10}$$

and the variance $\mathbb{V}[l_t]$ of $l_t$ as

$$\mathbb{V}[l_t] = \begin{cases} O(t^{3-\alpha}) & (1 < \alpha < 2) \\ O(t \ln t) & (\alpha = 2) \\ O(t) & (\alpha > 2) \end{cases} \quad \text{for } t >> 1, \tag{11}$$

suggesting the occurrence of a superdiffusion phenomenon if the flow velocity has long memory, i.e., if $\alpha \in (1,2)$ (see **Section A.1** in **Appendix**). This condition is satisfied in our application study in **Section 4**. This superdiffusion concerns the ensemble over the paths of velocity; for each fixed path, the eDNA profile is translated without spreading.

Discrete-time random flow velocity has also been studied in the context of traffic modeling (Böhme et al. 2025), and advective Lagrangian particle transport in random and possibly long-memory flow fields has been reported from an SPDE standpoint (Bevia et al. 2025; Cifani et al. 2026; Dorini and Cunha 2011); however, these studies have not been applied to more complex phenomena such as eDNA dynamics.

***Remark 3.*** We model flow velocity as a supOU process flow motivated by our application study, while our theoretical results hold true for any other càdlàg $u_t$ that is a.s. positive and bounded.

### 2.3 Mathematical analysis

We present a series of mathematical analysis results for our superprocess model. Throughout the rest of **Section 2** and the whole section of **Section 3**, we assume that a sample path of flow velocity $(u_t)_{0\le t\le T}$ is already given and is therefore formally seen as a nonnegative function on $[0,T]$, i.e., we employ the quenching approach (Böhme et al. 2025; Ji and Xiong 2022), where a path of $u$ is fixed.

In the remainder of this paper, we consider eDNA dynamics in a fixed time interval $[0,T]$ and suppose that **Assumption 1** presented below is satisfied.

***Assumption 1***

*The coefficients $R_t$ and $\sigma_t$ are continuous for all $t\in[0,T]$. Moreover, $f_t(x)$ is continuous, bounded for all $t\in[0,T]$ and $x\in\mathbb{R}$, and has compact support in $x$ with the length $\omega>0$. Also, $X_0$ has a compact support.*

This assumption is considered not restrictive in applications because flow velocity and decay should not be infinitely large. This property is satisfied for the supOU process. Considering an infinitely large noise intensity would pose technical difficulties in defining the noise term in the proposed model. The compact-support property of $f$ and $X_0$ is to avoid obtaining an infinite mass of eDNA in the river.

To continue, we define solutions to the superprocess model, which are defined through a martingale problem (e.g., Definition 1.20 in Etheridge (2000)).

***Definition 1***

*A nonnegative measure $X\in\mathbb{D}([0,T],\mathcal{M}_+)$ is said to be a weak solution to the superprocess model if, for any $\phi\in C_b^\infty$, the process $(M_t(\phi))_{0\le t\le T}$ defined by*

$$\begin{aligned} M_t(\phi) \equiv & \int_{\mathbb{R}}\phi(x)X_t(\mathrm{d}x)-\int_{\mathbb{R}}\phi(x)X_0(\mathrm{d}x) \\ & +\int_0^t\int_{\mathbb{R}}u_s\frac{\mathrm{d}\phi(x)}{\mathrm{d}x}X_s(\mathrm{d}x)\mathrm{d}s+\int_0^t\int_{\mathbb{R}}R_s\phi(x)X_s(\mathrm{d}x)\mathrm{d}s-\int_0^t\int_{\mathbb{R}}\phi(x)f_s(x)\mathrm{d}x\mathrm{d}s \end{aligned} \tag{12}$$

*is a martingale with quadratic variation*

$$\langle M_t(\phi)\rangle=\int_0^t\int_{\mathbb{R}}\sigma_s^2(\phi(x))^2X_s(\mathrm{d}x)\mathrm{d}s. \tag{13}$$

Our first theoretical result is about Laplace functionals, stating that the Laplace functional of the superprocess model must have a certain form and is given analytically.

***Proposition 1***

*For any nonnegative* $\psi \in C_b^\infty$ *and a weak solution* $X$ *to the superprocess model, it follows that*

$$\mathbb{E}\left[\exp\left(-\int_{\mathbb{R}}\psi(x)X_T(\mathrm{d}x)\right)\middle|\mathcal{F}_t\right]=\exp\left(-\mu_t-\int_{\mathbb{R}}\lambda_t(x)X_t(\mathrm{d}x)\right),\quad 0\le t\le T, \tag{14}$$

*where*

$$\lambda_t(x)=\frac{\psi\left(x-\int_t^T u_s\mathrm{d}s\right)\exp\left(-\int_t^T R_s\mathrm{d}s\right)}{1+\frac{1}{2}\psi\left(x-\int_t^T u_s\mathrm{d}s\right)\int_t^T\sigma_s^2\exp\left(-\int_s^T R_m\mathrm{d}m\right)\mathrm{d}s},\quad 0\le t\le T,\quad x\in\mathbb{R} \tag{15}$$

*and*

$$\mu_t=\int_t^T\left(\int_{\mathbb{R}}f_s(x)\lambda_s(x)\mathrm{d}x\right)\mathrm{d}s,\quad 0\le t\le T,\quad x\in\mathbb{R}. \tag{16}$$

Our second theoretical result concerns the well-posedness of the superprocess model, indicating that weak solutions are unique in terms of distribution.

***Proposition 2***

*There exists at least one weak solution to the superprocess model. Moreover, any its weak solution admits the common Laplace functional given in (14). Namely, the weak solution is unique in law.*

***Remark 4.*** Inspecting **Proofs of Propositions 1-2** suggests that, for the Laplace functional, "nonnegative $\psi \in C_b^\infty$" can be replaced by "sufficiently large $\psi \in C_b^\infty$ such that $\psi > -2\inf_{0\le t\le T} R_t / \sup_{0\le t\le T}\sigma_t^2$". Note that, for $\psi < 0$, the right-hand side of (14) is no longer bounded; this extension is rigorously justified under (38) and (47) by **Proposition 4(b)**, as explained in **Remark 5**. These findings are key to studying model uncertainties in **Section 3**.

## 3. Model uncertainty

We introduce model uncertainties in our superprocess model in the context of Hansen and Sargent (2001). The problem considered here is to estimate the (weighted) average of the eDNA concentration when the models are not necessarily correct. This problem would be encountered, for example, when estimating the total amount or average of eDNA concentration in a domain. More specifically, model uncertainties in the context of Hansen and Sargent (2001) refer to the following situation, where the SPDE (1) is a model identified from data, while the true model has an added drift, called distortion $d_t(x)$ here:

$$\frac{\partial X_t}{\partial t}=u_t\frac{\partial X_t(x)}{\partial x}-R_tX_t(x)+f_t(x)+\underbrace{d_t(x)}_{\text{Distortion}}+\sigma_t\sqrt{X_t(x)}\dot{W}. \tag{17}$$

We want to find the worst-case distortions $d_t(x)$ that give the maximum and minimum of an objective function, such as the eDNA concentration at some time instance, subject to the error size measured through the relative entropy between the benchmark and distorted models. The worst-case model turns out to be again a superprocess, and hence the tractability of the model is not degraded. The evaluation of model uncertainties in this framework has been mainly applied to SDEs (e.g., Guan and Liang 2026; Hansen and Souganidis 2025; Yang and Sun 2025), and to SPDEs only with additive noise (Baltas et al. 2019).

### 3.1 Girsanov transform

We consider Dawson's Girsanov transform (p. 135 in Etheridge (2000)):

$$Z_t = \exp\left( \int_0^t \int_{\mathbb{R}} \eta_s(X_s, x) Q(\mathrm{d}s, \mathrm{d}x) - \frac{1}{2} \int_0^t \int_{\mathbb{R}} \int_{\mathbb{R}} \eta_s(X_s, x) \eta_s(X_s, y) S(X_s; \mathrm{d}x, \mathrm{d}y) \mathrm{d}s \right), \quad 0 < t \le T, \tag{18}$$

where $\eta$ is a mapping from $[0,T] \times \mathcal{M}_+ \times \mathbb{R} \to \mathbb{R}$ that is adapted to the filtration generated by $X$ and $Q(\mathrm{d}s, \mathrm{d}x)$ is an orthogonal martingale measure (p.134 in Etheridge (2000)) with the covariance structure

$$Q(\mathrm{d}t, \mathrm{d}x) Q(\mathrm{d}t, \mathrm{d}y) = \delta(x - y) X_t(\mathrm{d}y) \mathrm{d}t \ \left(\equiv S(X_t; \mathrm{d}x, \mathrm{d}y) \mathrm{d}t\right), \quad 0 < t \le T \tag{19}$$

where $\delta$ is Dirac's Delta. Afterward, we define a new probability measure $\mathbb{Q}$ such that

$$\left.\frac{\mathrm{d}\mathbb{Q}}{\mathrm{d}\mathbb{P}}\right|_{\mathcal{F}_t} = Z_t, \quad 0 < t \le T. \tag{20}$$

According to pp. 135-136 in Etheridge (2000), if $\eta$ is strictly bounded and measurable, then $Z$ is a martingale under $\mathbb{P}$, $\mathbb{Q}$ is equivalent to $\mathbb{P}$, and $M_{\mathbb{Q},t}(\phi)$ ($0 \le t \le T$) given by

$$\begin{aligned} M_{\mathbb{Q},t}(\phi) \equiv & \int_{\mathbb{R}} \phi(x) X_t(\mathrm{d}x) - \int_{\mathbb{R}} \phi(x) X_0(\mathrm{d}x) - \underbrace{\int_0^t \int_{\mathbb{R}} \phi(x) \sigma_s \eta_s(X_s, x) X_s(\mathrm{d}x) \mathrm{d}s}_{\text{Distortion}} \\ & + \int_0^t \int_{\mathbb{R}} u_s \frac{\mathrm{d}\phi(x)}{\mathrm{d}x} X_s(\mathrm{d}x) \mathrm{d}s + \int_0^t \int_{\mathbb{R}} R_s \phi(x) X_s(\mathrm{d}x) \mathrm{d}s - \int_0^t \int_{\mathbb{R}} \phi(x) f_s(x) \mathrm{d}x \mathrm{d}s \end{aligned} \tag{21}$$

is a martingale for any $\phi \in C_b^\infty$ under $\mathbb{Q}$. The expectation under $\mathbb{Q}$ is denoted by $\mathbb{E}_{\mathbb{Q}}$. The quadratic variation of $M_{\mathbb{Q},t}(\phi)$ is given by (13).

### 3.2 Relative entropy and objective functions

Following the finite-dimensional case (e.g., Hess (2019)), the relative entropy between $\mathbb{Q}$ and $\mathbb{P}$ based on $\mathbb{P}$ at time $t$ is given by

$$\begin{aligned} \mathrm{REnt}_t(\mathbb{P}, \mathbb{Q}) &\equiv \mathbb{E}_{\mathbb{Q}}\left[ \ln\left( \left.\frac{\mathrm{d}\mathbb{Q}}{\mathrm{d}\mathbb{P}}\right|_{\mathcal{F}_t} \right) \right] \\ &= \mathbb{E}_{\mathbb{Q}}[\ln Z_t] \\ &= \frac{1}{2} \mathbb{E}_{\mathbb{Q}}\left[ \int_0^t \int_{\mathbb{R}} \int_{\mathbb{R}} \eta_s(X_s, x) \eta_s(X_s, y) S(X_s; \mathrm{d}x, \mathrm{d}y) \mathrm{d}s \right] \\ &= \frac{1}{2} \mathbb{E}_{\mathbb{Q}}\left[ \int_0^t \int_{\mathbb{R}} (\eta_s(X_s, x))^2 X_s(\mathrm{d}x) \mathrm{d}s \right]. \end{aligned} \tag{22}$$

With this formulation of relative entropy, based on the stochastic control framework where distortions are considered control variables to be optimized by nature (e.g., Hansen and Sargent 2001), we consider the following lower- and upper-bound problems of weighted eDNA concentration at time $T$ ; for fixed $\psi \in C_b^\infty$ ,

$$\textbf{(Lower-bound case)} \qquad \underline{\Phi} = \inf_{\eta} \left\{ \underbrace{\mathbb{E}_{\mathbb{Q}}\left[\int_{\mathbb{R}} \psi(x) X_T(\mathrm{d}x)\right]}_{\text{Expectation}} + \underbrace{w\mathrm{REnt}_T(\mathbb{P},\mathbb{Q})}_{\text{Entropic penalty}} \right\} \tag{23}$$

and

$$\textbf{(Upper-bound case)} \qquad \bar{\Phi} = \sup_{\eta} \left\{ \underbrace{\mathbb{E}_{\mathbb{Q}}\left[\int_{\mathbb{R}} \psi(x) X_T(\mathrm{d}x)\right]}_{\text{Expectation}} - \underbrace{w\mathrm{REnt}_T(\mathbb{P},\mathbb{Q})}_{\text{Entropic penalty}} \right\}, \tag{24}$$

where $w > 0$ is a weighting factor that balances the two terms in each objective function and $\eta$ in the "sup" and "inf" are explored for bounded and measurable mappings from $[0,T] \times \mathcal{M}_+ \times \mathbb{R} \to \mathbb{R}$ that are adapted to the filtration generated by $X$ . The function $\psi$ can be chosen so that it has a large value on a priority measurement area of the eDNA of aquatic species.

The problems (24) and (23) are infinite-dimensional stochastic control problems that may be difficult to address directly. We therefore employ the strategy to first consider and solve their finite-dimensional versions and guess their limits, as shown below.

### 3.3 Formal results

We consider the semidiscrete system with distortion:

$$\mathrm{d}X_t^{(i)} = \left( u_t \frac{X_t^{(i+1)} - X_t^{(i)}}{h} - R_t X_t^{(i)} + f_t^{(i)} + \sigma_t X_t^{(i)} \eta_t^{(i)} \right) \mathrm{d}t + \sigma_t \sqrt{X_t^{(i)}} \frac{\mathrm{d}\tilde{B}_t^{(i)}}{\sqrt{h}}, \quad t > 0, \quad i \in \mathbb{Z} \tag{25}$$

subject to the initial condition $X_0^{(i)} \geq 0$ ($i \in \mathbb{Z}$), where $\left(\tilde{B}_t^{(i)}\right)_{t \geq 0}$ is a standard one-dimensional Brownian motion on $\mathbb{Q}$ and all the elements of $\left\{\tilde{B}_t^{(i)}\right\}_{i \in \mathbb{Z}}$ are mutually independent.

For the lower-bound case, in view of the dynamic programming principle (Chapter 5.1 in Øksendal and Sulem (2019)), we infer that $\underline{\Phi}$ is approximated by $\varphi(0,\{x_i\})$, where $\varphi = \varphi(t,\{x_i\})$ solves the dynamic programming equation:

$$\frac{\partial \varphi}{\partial t} + \inf_{\{\eta_i\}} \left\{ \begin{array}{l} \sum_{i \in \mathbb{Z}} \left( u_t \frac{x_{i+1} - x_i}{h} - R_t x_i + f_t^{(i)} + \sigma_t x_i \eta_i \right) \frac{\partial \varphi}{\partial x_i} \\ + \frac{1}{2h} \sigma_t^2 \sum_{i \in \mathbb{Z}} x_i \frac{\partial^2 \varphi}{\partial x_i^2} + \frac{1}{2} w \sum_{i \in \mathbb{Z}} x_i \eta_i^2 h \end{array} \right\} = 0, \quad 0 \leq t < T, \quad x_i > 0 \quad (i \in \mathbb{Z}) \tag{26}$$

subject to the terminal condition $\varphi(T,\{x_i\}) = \psi_i x_i$ ($i \in \mathbb{Z}$) with $\eta_i \in \mathbb{R}$ ($i \in \mathbb{Z}$). We have

$$\eta_i^* = \arg\min_{\{\eta_i\}} \left\{ \sum_{i\in\mathbb{Z}} \sigma_t x_i \eta_i \frac{\partial \varphi}{\partial x_i} + \frac{1}{2} w \sum_{i\in\mathbb{Z}} x_i \eta_i^2 h \right\} = -\frac{\sigma_t}{wh} \frac{\partial \varphi}{\partial x_i}, \quad (27)$$

reducing (26) to

$$\frac{\partial \varphi}{\partial t} + \sum_{i\in\mathbb{Z}} \left( u_t \frac{x_{i+1} - x_i}{h} - R_t x_i + f_t^{(i)} \right) \frac{\partial \varphi}{\partial x_i} + \frac{1}{2h} \sigma_t^2 \sum_{i\in\mathbb{Z}} x_i \frac{\partial^2 \varphi}{\partial x_i^2} - \frac{\sigma_t^2}{2wh} \sum_{i\in\mathbb{Z}} x_i \left( \frac{\partial \varphi}{\partial x_i} \right)^2 = 0 . \quad (28)$$

Next, we guess the solution of the form with time-dependent coefficients $\tilde{F}$ and $G^{(i)}$:

$$\varphi(t, \{x_i\}) = \tilde{F}_t + \sum_{i\in\mathbb{Z}} x_i G_t^{(i)} h \quad \text{and hence} \quad \frac{\partial \varphi}{\partial x_i} = G_t^{(i)} h , \quad \frac{\partial^2 \varphi}{\partial x_i^2} = 0 . \quad (29)$$

We have $\eta_i^* = -\sigma_t G_t^{(i)} / w$. Substituting (29) into (28) yields

$$\frac{\mathrm{d}\tilde{F}_t}{\mathrm{d}t} + \sum_{i\in\mathbb{Z}} x_i \frac{\mathrm{d}G_t^{(i)}}{\mathrm{d}t} h + \sum_{i\in\mathbb{Z}} \left( u_t \frac{x_{i+1} - x_i}{h} - R_t x_i + f_t^{(i)} \right) G_t^{(i)} h - \frac{\sigma_t^2}{2w} \sum_{i\in\mathbb{Z}} x_i \left( G_t^{(i)} \right)^2 h = 0 , \quad (30)$$

leading to the semidiscrete system

$$\frac{\mathrm{d}G_t^{(i)}}{\mathrm{d}t} - u_t \frac{G_t^{(i)} - G_t^{(i-1)}}{h} - R_t G_t^{(i)} - \frac{\sigma_t^2}{2w} \left( G_t^{(i)} \right)^2 = 0 , \quad 0 \le t < T , \ (i \in \mathbb{Z}) \quad (31)$$

and

$$\frac{\mathrm{d}\tilde{F}_t}{\mathrm{d}t} + \sum_{i\in\mathbb{Z}} f_t^{(i)} G_t^{(i)} h = 0 , \quad 0 \le t < T \quad (32)$$

subject to the terminal conditions $\tilde{F}_T = 0$ and $G_T^{(\cdot)} = \psi_{\cdot}$.

This implies that under $h \to +0$, we formally obtain

$$\underline{\Phi} = F_0 + \int_{\mathbb{R}} G_0(x) X_0(\mathrm{d}x) , \quad (33)$$

where

$$G_t(x) = \frac{\psi\left( x - \int_t^T u_s \mathrm{d}s \right) \exp\left( -\int_t^T R_s \mathrm{d}s \right)}{1 + \frac{1}{2} \psi\left( x - \int_t^T u_s \mathrm{d}s \right) \int_t^T \frac{\sigma_s^2}{w} \exp\left( -\int_s^T R_m \mathrm{d}m \right) \mathrm{d}s} , \quad 0 \le t \le T , \quad x \in \mathbb{R} \quad (34)$$

and

$$F_t = \int_t^T \left( \int_{\mathbb{R}} f_s(x) G_s(x) \mathrm{d}x \right) \mathrm{d}s , \quad 0 \le t \le T . \quad (35)$$

Finally, the worst-case distortion leading to $\underline{\Phi}$ is obtained as $\eta_t^*(X_t, x) = -\sigma_t G_t(x) / w$, which is independent of $X$ and strictly bounded, showing that the worst-case distorted superprocess for the lower-bound case is the SPDE(1), with the decay rate $R_t$ being replaced by the larger one $R_t + \sigma_t^2 G_t(x) / w$. Here, we have

$$\overline{\Phi} = F_0 + \int_{\mathbb{R}} G_0(x) X_0(\mathrm{d}x) \quad (36)$$

with (35) and

$$G_t(x)=\frac{\psi\left(x-\int_t^T u_s\mathrm{d}s\right)\exp\left(-\int_t^T R_s\mathrm{d}s\right)}{1-\frac{1}{2}\psi\left(x-\int_t^T u_s\mathrm{d}s\right)\int_t^T\frac{\sigma_s^2}{w}\exp\left(-\int_s^T R_m\mathrm{d}m\right)\mathrm{d}s},\quad 0\le t\le T,\quad x\in\mathbb{R}. \tag{37}$$

By **Remark 4**, we can also address the upper-bound case in a similar way where we need to choose a sufficiently large $w$ so that

$$\psi(x)<2w\left(\inf_{0\le t\le T}R_t/\sup_{0\le t\le T}\sigma_t^2\right),\quad x\in\mathbb{R}. \tag{38}$$

If (38) is satisfied, then the worst-case distortion leading to $\bar{\Phi}$ is $\eta_t^*(X_t,x)=\sigma_t G_t(x)/w$, which is independent of $X$ and strictly bounded, showing that the worst-case distorted superprocess for the upper-bound case is formally the SPDE (1) with the decay rate $R_t$ being replaced by the smaller one $R_t-\sigma_t^2 G_t(x)/w$. The condition (38) implies that the aversion against uncertainties should be sufficiently small ( $w$ is large); without this condition, we encounter an unbounded (and hence not meaningful) $\bar{\Phi}$.

### 3.4 Verification

By the heuristic analysis, for the lower-bound case, we guess

$$\Phi(t,g)=F_t+\int_{\mathbb{R}}G_t(x)g(\mathrm{d}x),\quad 0\le t<T,\quad g\in\mathcal{M}_+. \tag{39}$$

In particular, (26) implies that the dynamic programming equation for the superprocess model is as follows (e.g., Section 5.2 in Ocello (2025)).

$$\frac{\partial\Phi}{\partial t}+\inf_{\eta}\left\{\begin{array}{l}\int_{\mathbb{R}}f_t(x)\delta\Phi\mathrm{d}x+\int_{\mathbb{R}}\left(-u_tD_x\Phi+\left(-R_t+\sigma_t\eta_t(x)\right)\delta\Phi\right)g(\mathrm{d}x)\\ +\frac{\sigma_t^2}{2}\int_{\mathbb{R}}\delta\Phi g(\mathrm{d}x)+\frac{1}{2}w\int_{\mathbb{R}}\left(\eta_t(x)\right)^2g(\mathrm{d}x)\end{array}\right\}=0,\quad 0\le t<T,\quad g\in\mathcal{M}_+ \tag{40}$$

subject to the terminal condition $\Phi(T,g)=\int_{\mathbb{R}}\psi(x)g(\mathrm{d}x)$ ( $g\in\mathcal{M}_+$ ), where we set the following derivatives of $\Phi$ (e.g., Definitions 4.3 and 4.4 in Ocello (2025)):

**(Partial derivative)**
$$\frac{\partial\Phi(t,g)}{\partial t}=\frac{\mathrm{d}F_t}{\mathrm{d}t}+\int_{\mathbb{R}}\frac{\partial G_t(x)}{\partial t}g(\mathrm{d}x), \tag{41}$$

**(Linear (flat) derivative 1)**
$$\delta\Phi(t,x,g)=G_t(x), \tag{42}$$

**(Linear (flat) derivative 2)**
$$\delta^2\Phi(t,x,g)=0, \tag{43}$$

and

**(Intrinsic derivative)**
$$D_x\Phi(t,x,g)=\partial_x\delta\Phi(t,x,g)=\frac{\partial G_t(x)}{\partial x}. \tag{44}$$

Afterward, we obtain

$$\eta_t^*(x)=-\frac{\sigma_t}{w}\delta\Phi(t,x,g),\quad 0\le t<T,\quad g\in\mathcal{M}_+. \tag{45}$$

Substituting (39) into (40) yields

$$\frac{\mathrm{d}F_t}{\mathrm{d}t} + \int_{\mathbb{R}} f_t(x) G_t(x) \mathrm{d}x + \int_{\mathbb{R}} \left( \frac{\partial G_t(x)}{\partial t} - u_t \frac{\partial G_t(x)}{\partial x} - R_t G_t(x) - \frac{\sigma_t}{2} \left( G_t(x) \right)^2 \right) g(\mathrm{d}x) = 0 , \quad (46)$$

which is solved by (34) and (35).

The guessed solution (39) is twice flat differentiable, and the partial differential (39) is defined a.e. because $u$ has bounded variation. The HJB equation is formally derived following Ocello (2025). Since our model, which contains decay, immigration, and controls acting on the mass, is not directly covered by the verification theorem therein, we verify the optimality directly in **Proof of Proposition 3**.

***Proposition 3***

***(a) Lower-bound case***

*Under **Assumption 1**, it follows that* $\underline{\Phi} = \Phi(0, X_0)$ *where* $G_t(x)$ *is given by (34) and* $\Phi$ *by (39). Moreover, (33) holds true and* $\eta_t^*(x) = -\frac{\sigma_t}{w} G_t(x)$ ($x \in \mathbb{R}$).

***(b) Upper-bound case***

*Under **Assumption 1** and (38), it follows that* $\bar{\Phi} = \Phi(0, X_0)$ *where* $G_t(x)$ *is given by (37) and* $\Phi$ *by (39). Moreover, (36) holds true and* $\eta_t^*(x) = \frac{\sigma_t}{w} G_t(x)$ ($x \in \mathbb{R}$).

We conclude this section by stating that $\eta^*$ defines a measure $\mathbb{Q}$ equivalent to $\mathbb{P}$.

***Proposition 4***

***(a) Lower-bound case***

*The probability measure* $\mathbb{Q}$ *induced by* $\eta = \eta^*$ *is equivalent to* $\mathbb{P}$.

***(b) Upper-bound case***

*The probability measure* $\mathbb{Q}$ *induced by* $\eta = \eta^*$ *is equivalent to* $\mathbb{P}$ *if (38) and*

$$R_t - \sigma_t \eta_t^* \geq 0, \quad 0 \leq t \leq T, \quad x \in \mathbb{R}. \quad (47)$$

***Remark 5.*** In view of **Propositions 3-4**, we derive the following formulae, which are versions of certainty equivalent risk measures (e.g., Vinel and Krokhmal 2017; Koski and Karvanen 2025):

$$\underline{\Phi} = -w \ln\left( \mathbb{E}_{\mathbb{P}} \left[ \exp\left( -\frac{1}{w} \int_{\mathbb{R}} \psi(x) X_T(\mathrm{d}x) \right) \right] \right) \quad \left( \leq \mathbb{E}_{\mathbb{P}} \left[ \int_{\mathbb{R}} \psi(x) X_T(\mathrm{d}x) \right] \right) \quad (48)$$

and

$$\bar{\Phi} = w \ln\left( \mathbb{E}_{\mathbb{P}} \left[ \exp\left( \frac{1}{w} \int_{\mathbb{R}} \psi(x) X_T(\mathrm{d}x) \right) \right] \right) \quad \left( \geq \mathbb{E}_{\mathbb{P}} \left[ \int_{\mathbb{R}} \psi(x) X_T(\mathrm{d}x) \right] \right), \quad (49)$$

where the latter follows under (38) and (47) from $\mathbb{E}_{\mathbb{P}}\left[Z_T\right]=1$ in **Proposition 4(b)**, because applying Itô's formula to $\Phi\left(t,X_t\right)/w$ yields $Z_T=\exp\left(\left(\int_{\mathbb{R}}\psi\left(x\right)X_T\left(\mathrm{d}x\right)-\bar{\Phi}\right)/w\right)$ for $\eta=\eta^*$.

## 4. Application study

We apply the proposed model to eDNA dynamics in a river environment.

### 4.1 Numerical methods

If one wants to simulate sample paths of the superprocess model, then Yoshioka (2026a), which is based on the semidiscrete system combined with an unconditionally stable time discretization that splits the discretization into advection and noise steps, can be utilized. The scheme is based on the numerical method developed by Abi Jaber (2025), which utilizes local exact solutions to the associated Riccati equation of an affine process so that no truncation is necessary while it outperforms the Euler–Maruyama method in efficiency. Moreover, splitting schemes are effective for preserving the nonnegativity of solutions to S(P)DEs with multiplicative noise (e.g., Djurdjevac et al. 2026; Kamrani and Hausenblas 2025), and our case is no exception. Because our focus is on the Laplace functional, we use this scheme only for generating demonstrative sample paths, and we apply a common backward Euler method to (35), where the integrand is available in a closed form. The supOU process is discretized by using the method of Yoshioka (2026d), which is based on the exact sampling scheme of Qu et al. (2021) for scalar OU processes with tempered-stable Lévy measures.

Note that the numerical method explained above is unnecessary if one's interest is on statistics but not sample paths of eDNA concentration. This is because the latter approach does not need to discretize the SPDE with a non-Lipschitz diffusion coefficient whose computation is more complicated and less efficient than that based on Laplace functionals. Moreover, we do not have any pathwise convergence proofs for eDNA concentrations.

### 4.2 Study site and data

The study site is the Mitoya River in the Hii River System, which is a Class A river system in Eastern Shimane Prefecture, Japan (**Fig. 1**). The authors have been sampling water at the Kyu-Shita Bridge (35°18′06.87″N, 132°53′23.43″E) near the confluence point of the mainstream Hii River and Mitoya River. The Mitoya River is a tributary of the Hii River system; as it has no dams, it is characterized by a river environment that is not subject to the same degree of anthropogenic control as the mainstream Hii River, which does have a dam.

The fish species we focus on here is the diadromous fish *P. altivelis*, which is a major inland fishery resource in Japan and is currently among the most commercially important fish species in the Hii River system. The life history of this fish is one year, and they migrate between the river and sea; in the spring, juvenile fish begin to migrate upstream from the sea into rivers; in the summer, they grow into adults; in the autumn, the adults spawn and die; and the hatched young swim downstream to spend the winter in the sea (see also Tsukamoto and Uchida (1992) for their detailed life history).

The sampling interval of water for eDNA observation is approximately one week, and the data used are April to October in 2025 and April to August in 2026. Water sampling was conducted only when the water was not turbid, and the samples were immediately placed in a cold environment after collection. Real-time PCR based on the primers of Yamanaka and Minamoto (2016) was used to quantify the eDNA concentration of each sample. The procedure followed the Environmental DNA Sampling and Experiment Manual Ver 3.0 (Chapter 5 in The eDNA Society (2024)). We used glass fiber filter paper (aperture size of 0.25 μm) and a DN-Sure Tissue Mini Kit (Integral) for DNA extraction and the KOD SYBR qPCR/RT set (TOYOBO). We also collected hydraulic data to model flow velocity and water temperature. We obtained hourly water level data (from September 1, 2016, to August 31, 2026) and hourly discharge at the Shin-Mitoya observation station (35°17′38″N, 132°52′49″E) from the Water Information System[1], which is the nearest public hydrological observation station from the Kyu-Shita Bridge.

For the water temperature, we collected 20-min interval data near the Shin-Mitoya observation station (U-20, Onset) in 2025 and a small part of 2026; in 2026, an accident occurred in which the equipment could not be recovered. To address this missing data issue, we use data obtained by the Hii River Fisheries Cooperative near the Shin-Mitoya observation station, where water temperatures are observed using a simple thermometer at fixed times (AM07:00 and noon). With respect to flow velocity, hourly data are used to model the supOU process; however, since the modeling of water temperature focuses on trends, this approach is unlikely problematic. **Fig. 2** plots the collected environmental data used in this study.

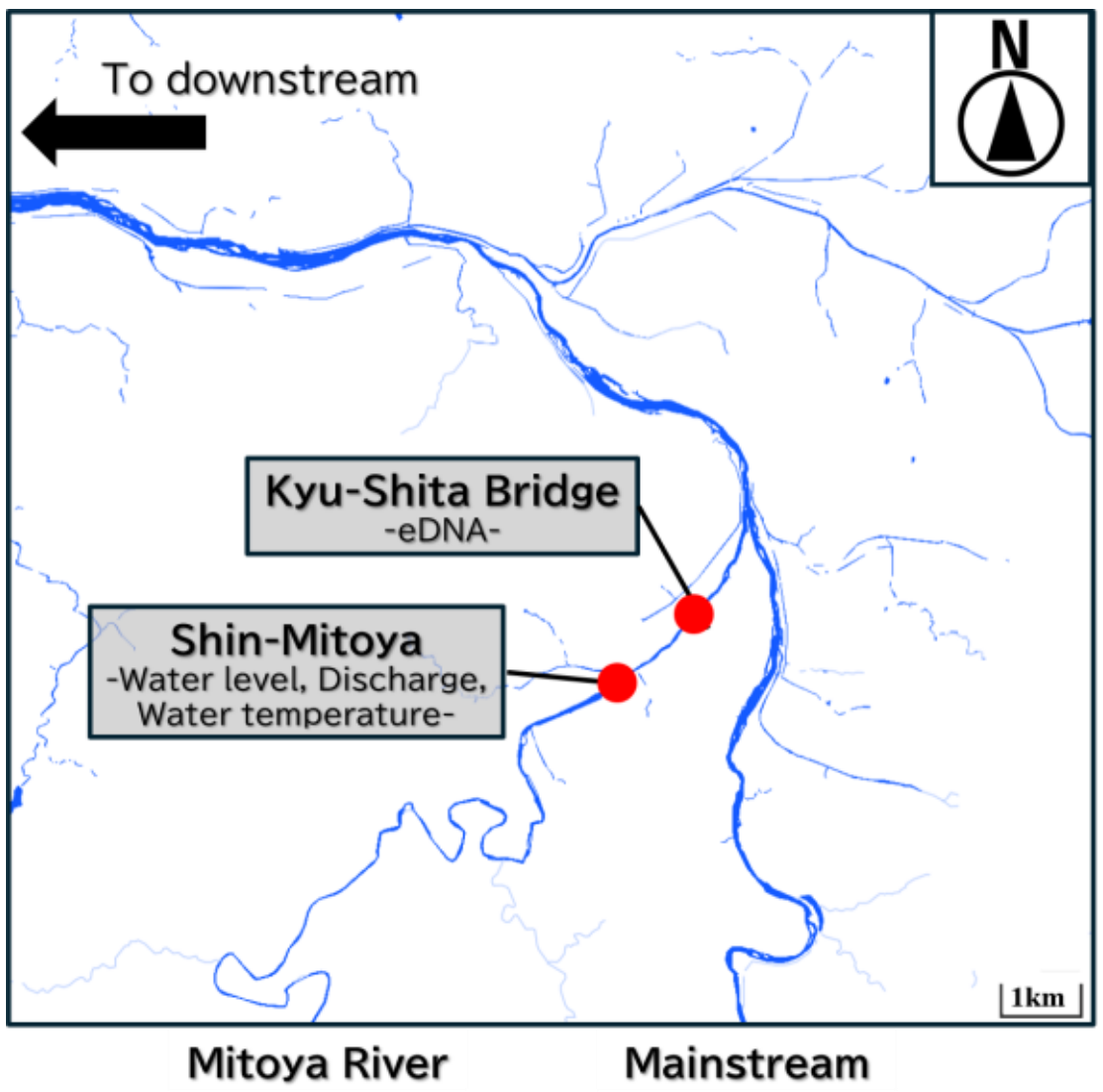

**Fig. 1** Map of the study area

[1] Water Information System. https://www1.river.go.jp/cgi-bin/SiteInfo.exe?ID=307041287705390 Last accessed on September 1 in 2026.

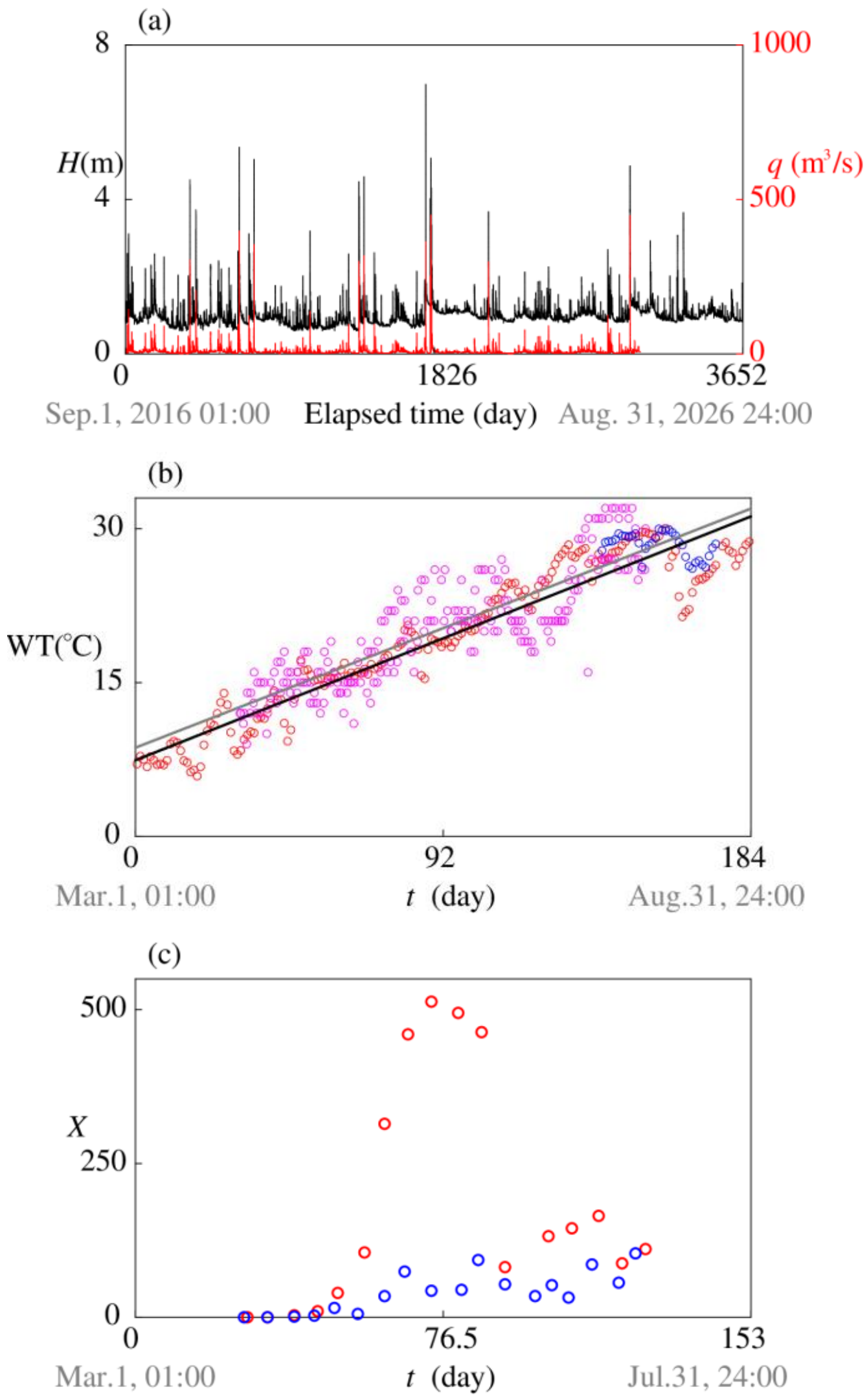

**Fig. 2** Collected data used in this study: (a) water level $H$ (m) and discharge $q$ (m$^3$/s) and (b) observed water temperature WT (°C) (circles) and their regressions (lines). For panel (b), we plot the daily averaged WT (°C) data in 2025 by U-20 (red), the daily averaged WT (°C) data in 2026 by U-20 (blue), and the instantaneous WT (°C) data by the Hii River Fisheries Cooperative (magenta). The black and gray lines are the best linear regressions for the WT (°C) data for 2025 and 2026, respectively. (c) eDNA concentration $X$ (copies/mL) at the Kyu-Shita Bridge from the beginning of March to the end of July in 2025 (red) and 2026 (blue)

### 4.3 Model setting

#### 4.3.1 Flow velocity

First, we estimate the hourly water velocity at the Shin-Mitoya observation station from the collected data as follows. We assume the following common Manning's formula, which is an empirical formula widely used in hydraulics and hydrology, for rectangular channels between flow velocity $u_t$ and water depth $h_t$ (Eq. (5.93) in Castro-Orgaz and Hager (2019)):

$$u_t = \frac{1}{n_{\mathrm{M}}} h_t^{\frac{2}{3}} I^{\frac{1}{2}}, \tag{50}$$

where $n_{\mathrm{M}} > 0$ (s/m$^{1/3}$) is the roughness coefficient and $I > 0$ (-) is the channel slope, both of which are unknown here. The discharge $q_t$ (m$^3$/s) is given by $u_t h_t b$, where $b > 0$ represents the river width and water level $H_t$ (m) is given by $h_t + \Delta$, with $\Delta$ (m) being a reference level. We rewrite (50) in terms of $q_t$ and $\eta_t$, both of which are available from the data:

$$q_t = \Psi \left( H_t - \Delta \right)^{\frac{5}{3}}, \tag{51}$$

where $\Psi = bI^{\frac{1}{2}} / n_{\mathrm{M}} > 0$ is considered one unknown constant in this paper. This equation has the two unknowns $\Delta > 0$ and $\Psi > 0$, which are identified by a nonlinear least-squares regression of (51) using the water level and discharge data from September 1 in 2016 to December 31 in 2024 under the constraint $H_t - \Delta \geq 0$ at all data points. Identified parameter values are $\Delta = 0.54$ (m) and $\Psi = 28.0$ (m$^{4/3}$/s). **Fig. 3** compares the observed and theoretical discharge-water depth relationships. **Fig. 4** plots the estimated hourly flow velocity.

Second, we identify a supOU process from data. The estimation procedure follows the least-square method against the autocorrelation and cumulants (average, variance, and skewness) employed in Yoshioka (2026b). The only difference is that they considered discharge while we consider flow velocity, both being nonnegative as time series data. We assume a gamma-type reversion measure $\pi(\mathrm{d}r) = p(r)\mathrm{d}r$ and a tempered-stable Lévy measure $v(\mathrm{d}z) = A_v e^{-b_v z} z^{-(\alpha_v + 1)} \mathrm{d}z$ ( $z > 0$ ), where $A_v > 0$ is the jump intensity parameter, $b_v > 0$ is the jump tilting parameter, and $\alpha_v < 1$ is the jump regularity parameter. The tempered stable model can capture both compound Poisson-type finite activity cases ($\alpha_v < 0$) and infinite activity cases ($0 \leq \alpha_v < 1$) in a unified way and is therefore considered a reasonable model for jump-driven long-memory processes (Yoshioka 2026b).

**Table 1** shows the estimated parameter values and **Table 2** compares empirical and theoretical cumulants. **Fig. 5** compares the fitted and theoretical autocorrelation functions and probability density functions (PDFs) of flow velocity. Particularly, the fitted model suggests that the flow velocity at the shin-Mitoya station has long memory ($1 < \alpha \leq 2$) and its jumps have infinite activities ($0 \leq \alpha_v < 1$) (**Fig. 5(b)**). This long-memory nature cannot be captured by models with exponential autocorrelation. For the PDF, both the empirical and fitted results have one mode at positive velocity, but they are slightly different in terms of the model height and width; the PDF of the fitted model is flatter than the empirical ones; nevertheless, the empirical cumulants are captured well by the model, where the average and variance are fitted with less than 1% errors and a skewness and kurtosis of 1–2% errors (**Table 2**).

Another key observation from the fitted model is that the long-memory nature ($1 < \alpha \leq 2$) of flow velocity along the Mitoya River. **Fig. 6** shows the theoretical variance of transported particles along the river, where its nonlinear behavior is clearly seen within hourly timescale.

**Table 1.** Estimated parameter values for flow velocity.

| Parameter | Estimated value |
|---|---|
| $\alpha$ (-) | 1.498 |
| $\beta$ (1/day) | 0.877 |
| $\alpha_v$ (m/s) | 0.785 |
| $b_v$ (s/m) | 0.660 |
| $A_v$ ($\mathrm{m}^{\alpha_v} / (\mathrm{s}^{\alpha_v} \cdot \mathrm{day})$) | 0.0411 |

**Table 2.** Comparison of empirical and theoretical cumulants.

| Parameter | Empirical | Theoretical | Relative error |
|---|---|---|---|
| Average (m/s) | 0.436 | 0.437 | 0.124% |
| Variance ($m^2/s^2$) | 0.0713 | 0.0712 | 0.0625% |
| Skewness (-) | 4.48 | 4.60 | 2.62% |
| Kurtosis (-)* | 43.9 | 43.3 | 1.30% |

* Kurtosis of Gaussian distributions is set to 0.

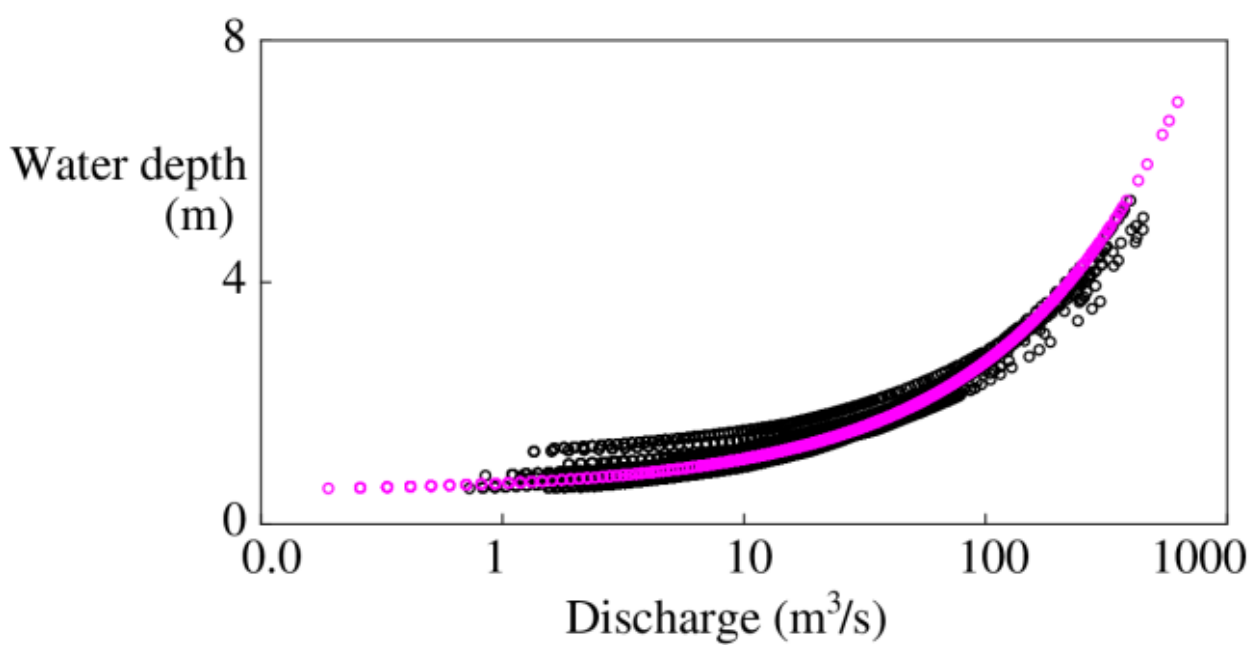

**Fig. 3** Comparison of observed (black) and theoretical (magenta) discharge–water depth relationships

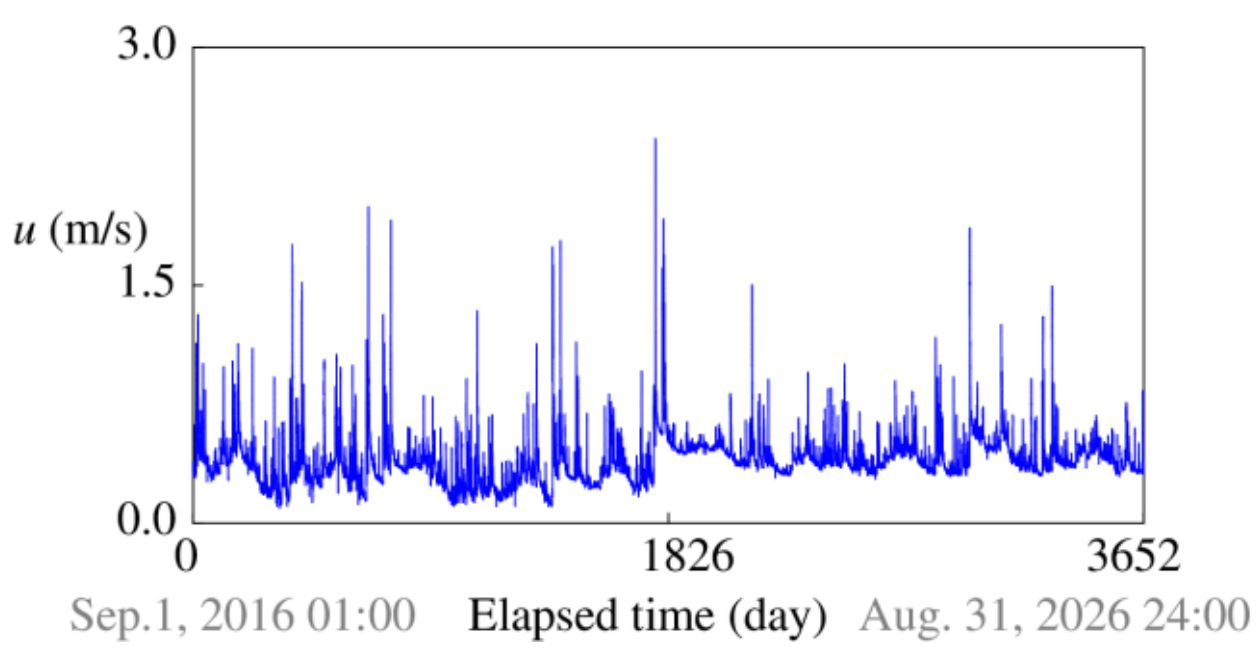

**Fig. 4** Estimated hourly flow velocity

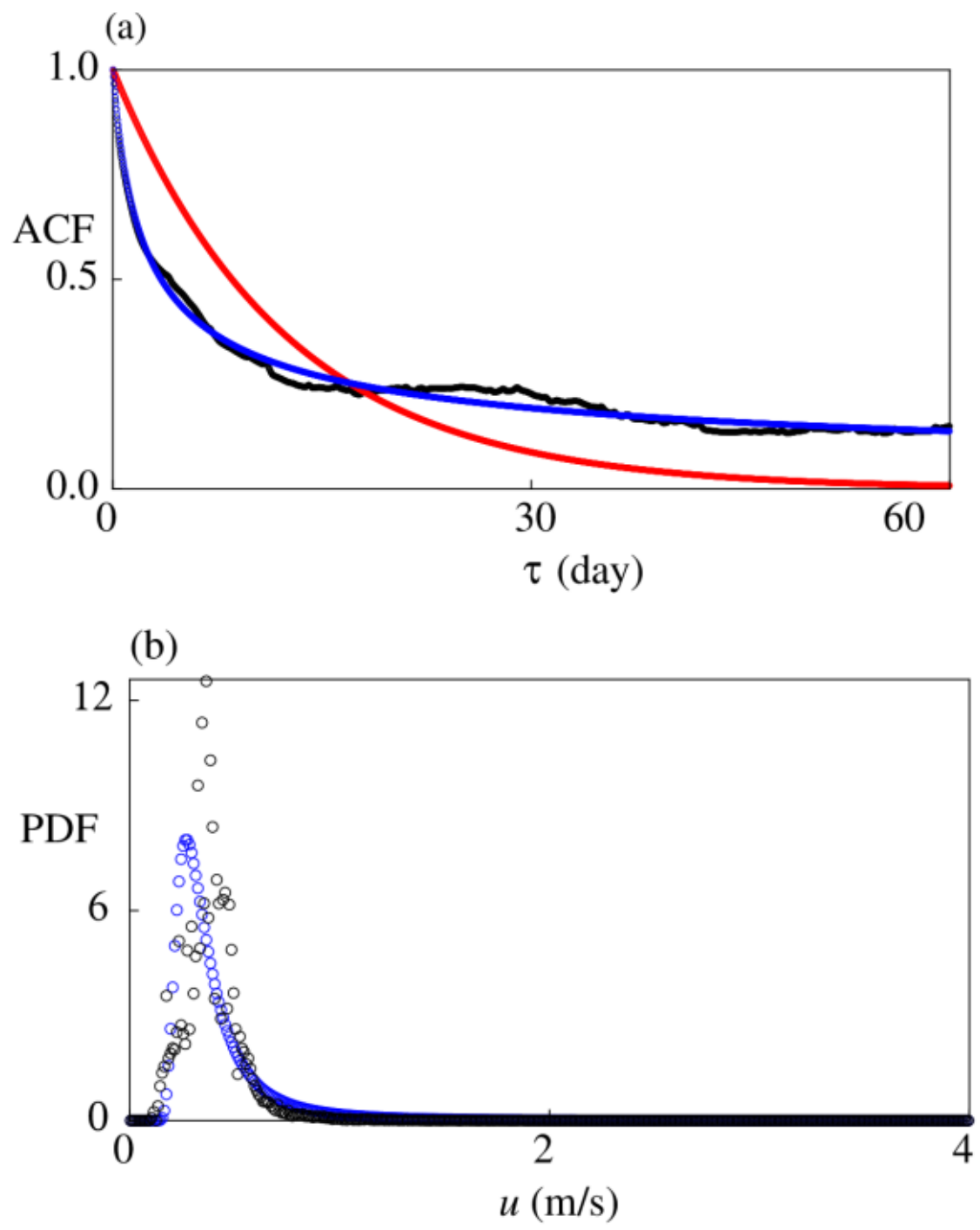

**Fig. 5** Comparisons of the empirical (black) and theoretical results (blue): (a) autocorrelation function of flow velocity where the exponential best fit $e^{-0.0812t}$ (red) is also plotted and (b) PDF of flow velocity. The theoretical PDF is computed by sampling the fitted model for 1,500,000 days

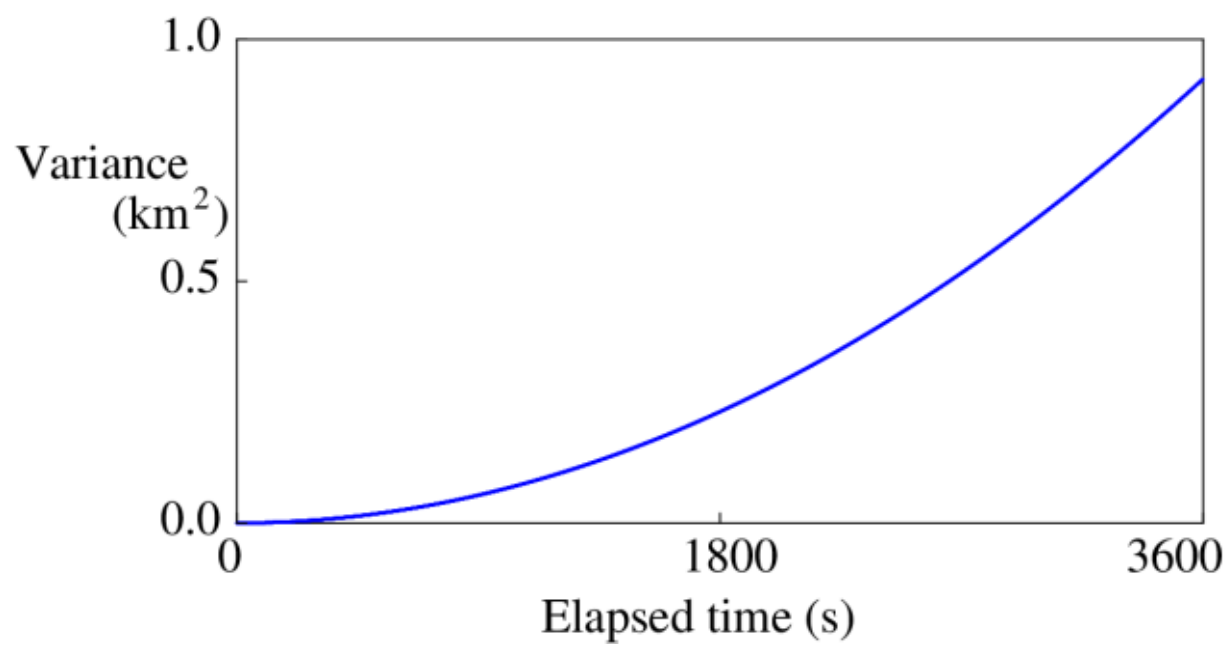

**Fig. 6** Theoretical variance for the transported length of particles along the Mitoya River

#### 4.3.2 Water temperature

We estimate the water temperature trend at the Shin-Mitoya observation station from the collected data as follows. Our focus is on the water temperature trend each year. Hence, the water temperature from March to August is modeled as follows:

$$\mathrm{WT}_t = c_0 + c_1 t , \tag{52}$$

where the time 0 here is the beginning of March, where $c_0, c_1$ are the parameters to be fitted. The linear trend model has been justified for daily-averaged water temperature data in the mainstream Hii River and another Class-A River in Yoshioka and Louriki (2026). The two parameters $c_0, c_1$ are fitted for each year and is plotted in **Fig. 2(b)**. **Table 3** summarizes the fitted parameter values. The results of the parameter estimation revealed that the slope of the WT trend is almost the same between 2025 and 2026, while the intercept is greater in 2026, resulting in a higher WT on average during spring to summer.

#### 4.3.3 Environmental DNA

For each year, the beginning of March 1 is set as $t = 0$. To fit the model against the eDNA data, we use the reconstructed flow velocity. For sensitivity and statistical analysis of the model, we use the identified supOU process as the flow velocity.

Environmental factors including water temperature have been suggested to affect the detection rate of eDNA (Tzafesta and Shokri 2025), and water type and water temperature (Lamb et al. 2022) as well as fish biomass (dos Anjos dos Santos et al. 2025) have been suggested to affect eDNA decay; these dependences have been regressed by exponential function models (Jo 2025). Considering data availability, for the decay rate $R_t$, following the simple experimental model of eDNA decay as a function of water temperature from 10 (°C) to 30 (°C) (Table 2 in Tsuji et al. (2017)), we set $R_t = 0.255\mathrm{WT}_t - 1.70$ (1/day). The positivity of this $R_t$ is guaranteed in our case. We assume a constant noise intensity $\sigma_t = \sigma = \text{constant}$ (copies$^{1/2}$km$^{1/2}$/(mL$^{1/2}$day$^{1/2}$)) unless otherwise specified, which does not affect the average eDNA concentration; in our numerical experiments, different values of $\sigma$ are examined. The initial condition is set as the zero measure $X_0 \equiv 0$ by assuming that the upstream migration of *P. altivelis* along the Mitoya River starts after March 1.

Finally, for source $f_t(x)$, we consider the following model:

$$f_t(x) = f_0 \max\left\{ \left( \frac{t - t_0 - x/V}{T_0} \right) \left( 1 - \frac{t - t_0 - x/V}{T_0} \right), 0 \right\}^{10} , \tag{53}$$

with parameters $f_0 > 0$ (source intensity), $t_0 \in \mathbb{R}$ (source shift), $V = 1$ (km/day) as the mean swimming velocity of individual *P. altivelis*, and $T_0 > 0$ (source duration) to be fitted; This model is a deterministic version of the source used in Yoshioka and Louriki (2026), where the power "10" comes from the average profile of daily fish count, which empirically matches observed eDNA data in the mainstream Hii River. This model assumes that the upstream migration of *P. altivelis* starts at the downstream end of the Mitoya River and ends at $t_0$ and $t_0 + T_0$, respectively. The constancy assumption of $V$ is due to the lack of data concerning the swimming of the fish in natural environments.

The empirical eDNA concentration data presented in **Fig. 2(c)** suggest a marked difference in eDNA concentrations between 2025 and 2026; specifically, while the peak of upstream migration is distinct in 2025 (around late May), it is indistinct in 2026. Furthermore, there is an approximately fivefold difference in the peak values. This implies that the migrating fish populations significantly differed between the two years. Nevertheless, it should be noted that such year-to-year variations, amounting to several-fold differences in the number of upstream-migrating *P. altivelis*, are not at all unnatural compared to cases in other rivers (e.g., Yoshioka 2026b). Because the focus in this paper is upstream migration, we fit the model against the data for 2025, with a clear peak in the eDNA concentration that corresponds to the upstream fish migration. In cases where there is no clear peak, as is the situation for 2026, a more complex modeling of the source term would be required; this will be a task for future work.

We fit the three parameters $f_0$, $t_0$, and $T_0$ so that the least-squares error between the theoretical average values of the eDNA concentration and the observed eDNA concentration values is minimized. Here, the theoretical average profile is used because we have only one empirical sample path of eDNA concentration each year. The influence of stochasticity is addressed in sensitivity analysis. We regard the observed data at the Kyu-Shita Bridge as a sample path $X_t(0)$. Because we do not have high-frequency or multiple sample paths, we estimate $f_0$, $t_0$, and $T_0$ by comparing the least-square errors between the observed eDNA concentration data and the average $\mathbb{E}\left[X_t(0)\right]$ as in Yoshioka (2026a), which in the present case is given by

$$\mathbb{E}\left[X_t(0)\right] = \int_0^t e^{-\int_s^t R_m \mathrm{d}m} f_s\left(\int_s^t u_m \mathrm{d}m\right)\mathrm{d}s, \tag{54}$$

where the expectation in (54) can be directly obtained by taking the expectation in (1). We examine the following parameter values to evaluate least-squares errors (also called root mean squares errors: RMSEs):

$$\left(f_0, t_0, T_0\right) = \left(10^6 \times (2100 + i),\ 10 + 0.1j,\ 115 + k\right) \quad \text{(copies/mL/day, day, day)}, \tag{55}$$

where $i = 1,2,3,\ldots,300$, $j = 1,2,3,\ldots,200$, and $k = 1,2,3,\ldots,10$. The range of parameter exploration has been determined through preliminary analysis with a wider exploration range under a coarser resolution. The duration to evaluate least-squares errors is from April 1 to the middle of June (June 9 in the present case) following Yoshioka (2026a). The fitted parameter values are shown in **Table 4**.

The fitted results are shown in **Fig. 8 in Appendix**, demonstrating that the fitted model captures the empirical results on average. According to the fitted model, the upstream migration of *P. altivelis* in 2025 along the Mitoya River started on March 26 and ended on July 22 in 2025. The least-squares error in this case is 65.4 (copies/mL), and this value should have the same order as the standard deviation of the eDNA concentration at the downstream end. We found that choosing $\sigma = 1$ (copies$^{1/2}$km$^{1/2}$/(mL$^{1/2}$day$^{1/2}$)) results in this condition and the use of this value unless otherwise specified; see also **Fig. 7(c)**.

#### 4.3.4 Computational conditions

To generate sample paths of our SPDE, we truncate the domain into the finite interval $(0, L)$ with $L = 36$ (km) based on an earlier study (Yoshioka 2026a), where the lower boundary $x = 0$ (km) corresponds to the Kyu-Shita Bridge, at which a free-outflow condition is specified. The no-flux condition of eDNA is specified at the upper boundary $x = L$ (km). The domain is divided into 1,200 grids. The time increment is 1/96 (day). For the demonstrative SPDE computation, we generate 10,000 sample paths. For the Laplace functionals of model uncertainties, we generate 100,000 sample paths.

**Table 3.** Fitted parameter values of the water temperature model. $R^2$ values are also presented.

| Parameter | 2025 | 2026 |
| --- | --- | --- |
| $c_0$ (°C) | 7.42 | 8.64 |
| $c_1$ (°C/day) | 0.129 | 0.127 |
| $R^2$ | 0.904 | 0.779 |

**Table 4.** Fitted parameter values of the source $f$ .

| Parameter | Value |
| --- | --- |
| $f_0$ (copies/mL/day) | $2360 \times 10^6$ |
| $t_0$ (day) | 25.8 |
| $T_0$ (day) | 118 |

### 4.4 Results and discussion

#### 4.4.1 Solution profile

For notational simplicity, **Fig. 7(a)** shows a computed sample path of the eDNA concentration under the assumption that $X$ admits a density, which is again denoted by $X$ here. **Figs. 7(b) and 7(c)** show the average (Ave) and standard deviation (Std) of $X$. **Fig. 9 in Appendix** shows a computed sample path of the flow velocity $u$. The results obtained suggest that the eDNA concentration evolves due to the space-time source term $f$ in a way that the front is moved toward upstream almost linearly with $u$.

The average of eDNA concentration has a profile analogous to that of the sample path but is seemingly smooth. The standard deviation of eDNA concentration also has a similar profile with a small peak near the upstream-end, which is considered caused by the finite-length domain. **Figs. 10 and 11 in Appendix** show the results with larger $\sigma$, demonstrating that sample paths become more irregular and intermittent with higher standard deviations for larger $\sigma$ values.

#### 4.4.2 Model uncertainty

We compute $\underline{\Phi}$ and $\overline{\Phi}$ for different settings of $\psi$ and $w$ to evaluate the influences of model uncertainties on statistics. We sample $\Phi$ obtained in **Proposition 3** for different paths of flow velocity $u$ because this proposition holds true for each fixed path of $u$. We fix $T = 90$ (day), which is around the middle of the migration period $[t_0, t_0 + T_0]$. We consider the case $\psi(x) = L^{-1}\mathbb{I}(x \in (0, L))$ such that the units of $\underline{\Phi}$ and $\overline{\Phi}$ become that of the eDNA concentration (copies/mL). This $\psi$ is not continuous, but we can choose a sequence $(\psi_\varepsilon)_{\varepsilon>0}$ of bounded and continuously differentiable functions such that regularizes $\psi_\varepsilon \underset{\varepsilon\to 0}{\to} \psi$ pointwise, and hence we can choose $\psi_\varepsilon$ that is arbitrary close to $\psi$. The coefficients $G_t(x)$ and $F_t$ are still well-defined with this $\psi$. For the upper-bound case, we need to set $w \geq 0.08$ to comply with the sufficient condition imposed in **Proposition 3**. The benchmark average value, which is defined as the average of benchmark limit $\lim_{w\to+\infty} \underline{\Phi}$ (or equivalently $\lim_{w\to+\infty} \overline{\Phi}$), is computed as 234.3 (copies/mL). The benchmark standard deviation value is computed as 46.1 (copies/mL).

**Tables 5 and 6** summarize the average, standard deviation, and coefficient of variation (standard deviation divided by average) for the lower- and upper-bound cases, respectively. The results obtained suggest that both the averages and the standard deviations of $\underline{\Phi}$ and $\overline{\Phi}$ depend continuously on $w$ and monotonically converge to the benchmark limit as $w \to +\infty$. The coefficient of variation of $\overline{\Phi}$ (resp., $\underline{\Phi}$) is decreasing in $w > 0$ (resp., increasing in $w > 0$). This implies that more uncertain models are more randomized (resp., less randomized) in a normalized sense for upper-bound cases (resp., lower-bound cases). The relative difference as a degree of relative difference between uncertain and benchmark models is also reported in **Tables 5 and 6**, yielding a rough estimate that the relative difference is $O(w^{-1})$ to $O(10w^{-1})$ % within the range of computational conditions examined here. From an application standpoint, based on Jo and Murakami (2026), measurement and data processing errors of eDNA variance for fish would be at most 10%, implying that realistic estimates of $w$ are $O\left(10^{-1}\right)$ to $O\left(10^{0}\right)$ under the present setting. However, these estimated values may vary depending on factors such as the target species and the range of eDNA concentration. Therefore, it is desirable to accumulate data across a diverse range of species and conditions. Nevertheless, our findings would provide guidelines on how to design $w$ under the assumption that model uncertainty can be discussed within the Hansen–Sargent framework.

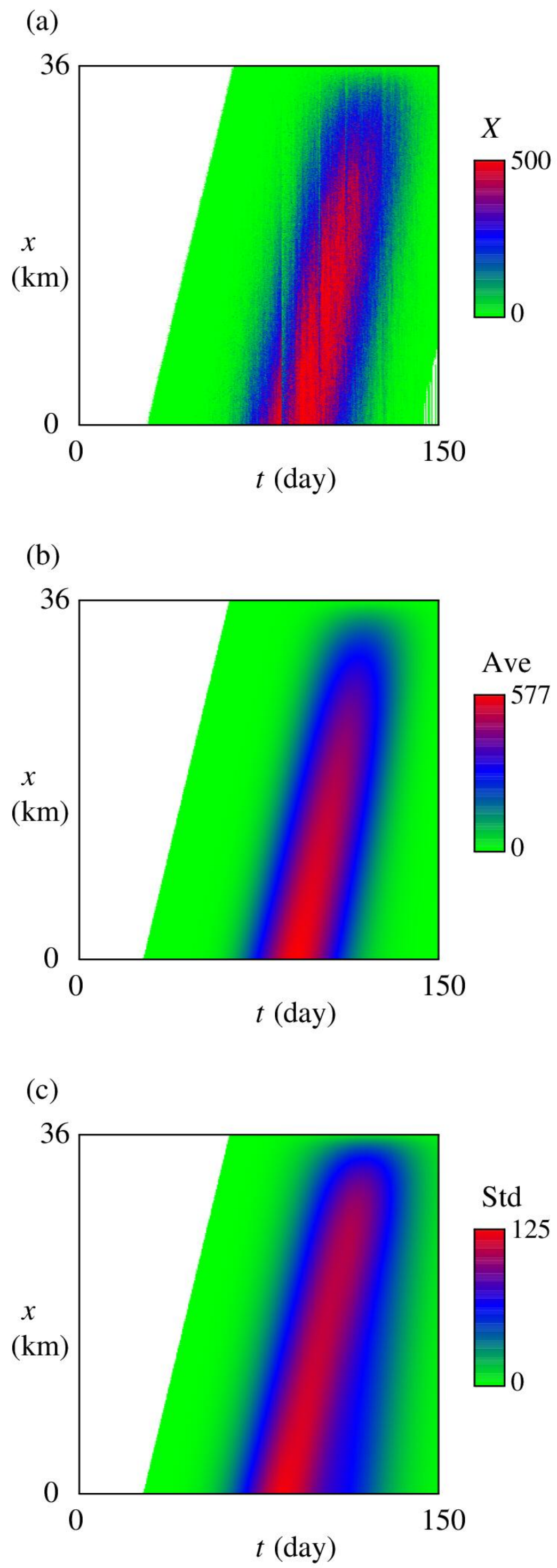

**Fig. 7** (a) Sample path of eDNA concentration $X$ (copies/mL), (b) average (Ave) of $X$, and (c) standard deviation (Std) of $X$ ($\max X = 1004$ (copies/mL) in the space-time domain). The white area represents zero value

**Table 5.** Average, standard deviation, coefficient of variation, and relative difference from average of $\underline{\Phi}$ from the benchmark case.

| $w$ | Average (copies/mL) | Standard deviation (copies/mL) | Coefficient of variation (-) | Relative Difference (%) |
|---|---|---|---|---|
| 0.5 | 211.0 | 39.6 | 0.188 | −9.92 |
| 1 | 221.8 | 42.8 | 0.193 | −5.32 |
| 2 | 227.8 | 44.6 | 0.196 | −2.76 |
| 4 | 231.0 | 45.5 | 0.197 | −1.41 |
| 8 | 232.6 | 45.9 | 0.197 | −0.71 |
| 16 | 233.4 | 46.3 | 0.198 | −0.36 |
| 32 | 233.8 | 46.4 | 0.198 | −0.18 |

**Table 6.** Same as **Table 5** but here, we present the results for $\bar{\Phi}$.

| $w$ | Average (copies/mL) | Standard deviation (copies/mL) | Coefficient of variation (-) | Relative Difference (%) |
|---|---|---|---|---|
| 0.5 | 266.6 | 57.0 | 0.214 | 13.81 |
| 1 | 248.9 | 51.2 | 0.206 | 6.26 |
| 2 | 241.3 | 48.7 | 0.202 | 2.99 |
| 4 | 237.7 | 47.6 | 0.200 | 1.47 |
| 8 | 236.0 | 47.1 | 0.199 | 0.73 |
| 16 | 235.1 | 46.8 | 0.199 | 0.36 |
| 32 | 234.7 | 46.7 | 0.199 | 0.18 |

## 5. Conclusions

We proposed a superprocess model for space-time eDNA dynamics in one-dimensional rivers and addressed its theoretical and computational analysis. The model assumes that the flow velocity and decay depend on the environment; in particular, the flow velocity is a long-memory process. The analytical tractability of the Laplace functional associated with the model played a role in the theoretical analysis of the model, ranging from uniqueness to uncertainty evaluation. Finally, the proposed model was applied to real eDNA data from the Hii River, Japan. Through this study, we deepened understanding of mathematical modeling based on superprocesses for eDNA dynamics and its applications, particularly in cases where the model is not necessarily correct.

This paper assumed that the flow velocity is time-dependent but not space-dependent for modeling simplicity, although this situation was found to be far from nontrivial. In such a case, the flow field can be considered governed by a hydrodynamic model driven by jump-driven input data, which is a difficult subject by itself. In our case, enough data for running the hydrodynamic model are not currently available for the Mitoya River, and we will address this issue in the future based on a joint measurement campaign of hydrodynamics and eDNA dynamics. Another topic to be addressed in the future is the generalization of the proposed model. For example, we assumed the noise to be space-time white, while one may consider generalizing the noise term so that the nonnegativity of solutions is guaranteed while more realistic stochasticity, such as colored space-time noise, can be considered. From an application standpoint, we focused on the upstream migration of *P. altivelis*, while modeling of the life history, including downstream migration, has not yet been achieved. We aim to overcome this challenge by continuing to collect eDNA data from the Hii River and other rivers. Regarding this, particularly the utilization of data acquired in 2026, the keys to success will likely be the formulation of more advanced models for fish school populations and the collection of further subsequent data.

The tractability of the proposed framework would carry over to SPDEs on metric graphs that approximate river networks (Lilleborge et al. 2026) if proper internal boundary conditions are specified at nodes. Using a time-fractional term (Modanli et al. 2023) or a Volterra formulation (Shen and Tsai 2024) is also a modeling option for representing memory effects, although the tractability of our model would not carry over to these cases.

Finally, the message we wish to convey through this study is that an industry-academia collaboration framework is vital in performing mathematical modeling for real-world problems. Such a collaboration would significantly help data collection, potentially leading to better mathematical models.

## Appendix

### A.1 Technical Proofs

***Proof of the estimate (11)***

In this proof, the expectation $\mathbb{E}$ means that for the randomness generated by the Poisson random measure $\mathfrak{N}$. By (7), we have

$$\begin{aligned}\mathbb{V}[l_t] &= \mathbb{E}\left[\left(\int_0^t (u_s - \mathbb{E}[u_s])\mathrm{d}s\right)^2\right] \\ &= \mathbb{E}\left[\int_0^t\int_0^t (u_s - \mathbb{E}[u_s])(u_m - \mathbb{E}[u_m])\mathrm{d}s\mathrm{d}m\right] \\ &= \mathbb{V}[u_0]\int_0^t\int_0^t \rho(|s-m|)\mathrm{d}s\mathrm{d}m \\ &= 2\mathbb{V}[u_0]\int_0^t\int_m^t \rho(s-m)\mathrm{d}s\mathrm{d}m \\ &= 2\mathbb{V}[u_0]\int_0^t\int_m^t (1+\beta(s-m))^{1-\alpha}\,\mathrm{d}s\mathrm{d}m.\end{aligned} \tag{A1}$$

If $\alpha \neq 2,3$, we have

$$\begin{aligned}\int_0^t\int_m^t (1+\beta(s-m))^{1-\alpha}\,\mathrm{d}s\mathrm{d}m &= \int_0^t \left[\frac{1}{(2-\alpha)\beta}(1+\beta(s-m))^{2-\alpha}\right]_m^t \mathrm{d}m \\ &= \frac{1}{(2-\alpha)\beta}\int_0^t \left((1+\beta(t-m))^{2-\alpha} - 1\right)\mathrm{d}m \\ &= \frac{1}{(2-\alpha)\beta}\left[-\frac{1}{(3-\alpha)\beta}(1+\beta(t-m))^{3-\alpha} - m\right]_0^t \\ &= \frac{1}{(2-\alpha)\beta}\left(\frac{(1+\beta t)^{3-\alpha} - 1}{(3-\alpha)\beta} - t\right).\end{aligned} \tag{A2}$$

For $t >> 1$, we obtain the asymptotic expansions

$$\mathbb{V}[l_t] \approx \frac{2\mathbb{V}[u_0]\beta^{1-\alpha}}{(2-\alpha)(3-\alpha)} t^{3-\alpha} \quad \text{if } 1 < \alpha < 2 \tag{A3}$$

and

$$\mathbb{V}[l_t] \approx \frac{2\mathbb{V}[u_0]}{(\alpha-2)\beta} t \quad \text{if } \alpha > 2\,. \tag{A4}$$

If $\alpha = 2,3$, then logarithmic terms appear in the estimates of $\mathbb{V}[l_t]$. If $\alpha = 2$, then

$$\int_0^t\int_m^t\left(1+\beta(s-m)\right)^{1-\alpha}\mathrm{d}s\mathrm{d}m=\int_0^t\int_m^t\left(1+\beta(s-m)\right)^{-1}\mathrm{d}s\mathrm{d}m$$
$$=\frac{1}{\beta}\int_0^t\ln\left(1+\beta(t-m)\right)\mathrm{d}m$$
$$=\frac{1}{\beta}\left[-\frac{1}{\beta}\left(1+\beta(t-m)\right)\ln\left(1+\beta(t-m)\right)-m\right]_0^t \tag{A5}$$
$$=\frac{1}{\beta}\left(\frac{1}{\beta}(1+\beta t)\ln(1+\beta t)-t\right).$$

For $t>>1$, we obtain

$$\mathbb{V}[l_t]\approx\frac{2\mathbb{V}[u_0]}{\beta}t\ln t \quad \text{if} \quad \alpha=2. \tag{A6}$$

If $\alpha=3$, then

$$\int_0^t\int_m^t\left(1+\beta(s-m)\right)^{1-\alpha}\mathrm{d}s\mathrm{d}m=\frac{1}{\beta}\int_0^t\left(1-\left(1+\beta(t-m)\right)^{-1}\right)\mathrm{d}m=\frac{1}{\beta}\left(t-\frac{1}{\beta}\ln(1+\beta t)\right), \tag{A7}$$

again yielding (A4).

□

***Proof of Proposition 1***

The proof is divided into several steps.

**Step 1. Guess governing system**

Here, we use the martingale problem in **Definition 1** with the time-dependent test function $\lambda_t$, which is smooth in $x$ and absolutely continuous in $t$ with a bounded derivative. This is justified by applying (12) on each subinterval of a partition of $[0,t]$ with the test function approximated by its time averages on a dyadic partition of $[0,t]$ (see, Proposition II.5.7 in Perkins (2002)); the time dependence of $u$ and $f$ does not affect this argument. Then, $\frac{\partial\lambda_s(x)}{\partial s}$ is added to the integrand in (12).

We express the right-hand side of (14) as $\exp(-U_t)$. Applying Itô's formula yields

$$\mathrm{d}\left(\exp(-U_t)\right)=-\exp(-U_t)\left(-\mathrm{d}U_t+\frac{1}{2}(\mathrm{d}U_t)^2\right), \tag{A8}$$

where

$$\begin{aligned}\mathrm{d}U_t&=\frac{\mathrm{d}\mu_t}{\mathrm{d}t}\mathrm{d}t+\int_{\mathbb{R}}\frac{\partial\lambda_t(x)}{\partial t}X_t(\mathrm{d}x)\mathrm{d}t+\int_{\mathbb{R}}\lambda_t(x)\mathrm{d}X_t(\mathrm{d}x)\\&=\frac{\mathrm{d}\mu_t}{\mathrm{d}t}\mathrm{d}t+\int_{\mathbb{R}}\frac{\partial\lambda_t(x)}{\partial t}X_t(\mathrm{d}x)\\&\quad-\int_{\mathbb{R}}u_t\frac{\partial\lambda_t(x)}{\partial x}X_t(\mathrm{d}x)\mathrm{d}t-\int_{\mathbb{R}}R_t\lambda_t(x)X_t(\mathrm{d}x)\mathrm{d}t+\int_{\mathbb{R}}\lambda_t(x)f_t(x)\mathrm{d}x+M_t(\lambda_t)\\&=\left(\frac{\mathrm{d}\mu_t}{\mathrm{d}t}+\int_{\mathbb{R}}\lambda_t(x)f_t(x)\mathrm{d}x+\int_{\mathbb{R}}\left(\frac{\partial\lambda_t(x)}{\partial t}-u_t\frac{\partial\lambda_t(x)}{\partial x}-R_t\lambda_t(x)\right)X_t(\mathrm{d}x)\right)\mathrm{d}t+M_t(\lambda_t)\end{aligned} \tag{A9}$$

and

$$\left(\mathrm{d}U_t\right)^2 = \mathrm{d}\left\langle M_t\left(\lambda_t\right)\right\rangle = \int_{\mathbb{R}} \sigma_t^2 \left(\lambda_t\left(x\right)\right)^2 X_t\left(\mathrm{d}x\right)\mathrm{d}t \ . \tag{A10}$$

In summary, we obtain

$$\begin{aligned} &\mathrm{d}\left(\exp\left(-U_t\right)\right) \\ &= -\exp\left(-U_t\right)\left(\left(\begin{aligned} &\frac{\partial \mu_t}{\partial t} + \int_{\mathbb{R}} \lambda_t\left(x\right) f_t\left(x\right)\mathrm{d}x \\ &+\int_{\mathbb{R}}\left(\frac{\partial \lambda_t\left(x\right)}{\partial t} - u_t \frac{\partial \lambda_t\left(x\right)}{\partial x} - R_t\lambda_t\left(x\right) - \frac{\sigma_t^2}{2}\left(\lambda_t\left(x\right)\right)^2\right) X_t\left(\mathrm{d}x\right) \end{aligned}\right)\mathrm{d}t + \mathrm{d}M_t\left(\lambda_t\right)\right). \end{aligned} \tag{A11}$$

The left-hand side of (14) becomes a martingale only if the drift term of (A11) vanishes, yielding a time-backward system of equations to find $\mu$ and $\lambda$ :

$$\frac{\partial \lambda_t\left(x\right)}{\partial t} - u_t \frac{\partial \lambda_t\left(x\right)}{\partial x} - R_t\lambda_t\left(x\right) - \frac{\sigma_t^2}{2}\left(\lambda_t\left(x\right)\right)^2 = 0 \ , \ \ 0 \le t < T \ , \ \ x \in \mathbb{R} \tag{A12}$$

and

$$\frac{\partial \mu_t}{\partial t} + \int_{\mathbb{R}} f_t\left(x\right)\lambda_t\left(x\right)\mathrm{d}x = 0 \ , \ \ 0 \le t < T \ , \tag{A13}$$

where the terminal conditions are $\mu_T = 0$ and $\lambda_T = \psi$ . We show that $\exp\left(-U_t\right)$ is truly a martingale due to $U_t \ge 0$ , which is obtained from nonnegativity of $\mu$ and $\lambda$ .

**Step 2. Solve the hyperbolic balance law**

We show that the (time-backward) hyperbolic balance law (A12) subject to the terminal condition $\lambda_T = \psi$ is uniquely solvable. First, by the transformation $\left(\theta, y\right) = \left(T - t, x - \int_t^T u_s \mathrm{d}s\right)$ (this transformation is allowed because $u_s$ is a.s. positive and bounded), by denoting $\rho_\theta\left(y\right) = \lambda_t\left(x\right)$ , we obtain

$$\frac{\partial \rho_\theta\left(y\right)}{\partial \theta} = -R_{T-\theta}\rho_\theta\left(y\right) - \frac{\sigma_{T-\theta}^2}{2}\left(\rho_\theta\left(y\right)\right)^2 , \ \ 0 < \theta \le \Theta_T , \ \ y \in \mathbb{R} \tag{A14}$$

subject to the initial condition $\rho_0 = \psi\left(y\right)$ ( $y \in \mathbb{R}$ ), where $\Theta_T = \int_0^T u_s \mathrm{d}s$ . From this ordinary differential equation (ODE), we obtain the *a priori* estimate

$$0 \le \rho_\theta\left(y\right) \le \psi\left(y\right), \ \ 0 \le \theta \le \Theta_T , \ \ y \in \mathbb{R} , \tag{A15}$$

which should be satisfied at each point for all reasonable solutions (e.g., entropy solutions) to (A12). Indeed, substituting $\rho_\theta\left(y\right) = 0$ into (A14) results in $\frac{\partial \rho_\theta\left(y\right)}{\partial \theta} = 0$ . We also have

$$\frac{\partial \rho_\theta\left(y\right)}{\partial \theta} \le -R_{T-\theta}\rho_\theta\left(y\right) \to \rho_\theta\left(y\right) \le \psi\left(y\right)\exp\left(-\int_0^\theta R_{T-m}\mathrm{d}m\right) \le \psi\left(y\right) , \ \ 0 < \theta \le \Theta_T , \ \ y \in \mathbb{R} . \tag{A16}$$

Moreover, (A14) is explicitly solved as follows:

$$\rho_\theta(y)=\frac{\psi(y)\exp\left(-\int_0^\theta R_{T-m}\mathrm{d}m\right)}{1+\frac{1}{2}\psi(y)\int_0^\theta \sigma_{T-m}^2\exp\left(-\int_0^m R_{T-\tau}\mathrm{d}\tau\right)\mathrm{d}m},\quad 0\le\theta<\Theta_T,\quad y\in\mathbb{R}. \tag{A17}$$

Indeed, applying the transformation $\tilde{\rho}_\theta(y)=1/\rho_\theta(y)$ yields a linear ODE of $\tilde{\rho}$ that can be solved exactly, and after solving this equation, we return to $\rho_\theta(y)$. Finally, returning to the original parameters $(t,x)$ yields

$$\lambda_t(x)=\frac{\psi\left(x-\int_t^T u_s\mathrm{d}s\right)\exp\left(-\int_t^T R_s\mathrm{d}s\right)}{1+\frac{1}{2}\psi\left(x-\int_t^T u_s\mathrm{d}s\right)\int_t^T \sigma_s^2\exp\left(-\int_s^T R_m\mathrm{d}m\right)\mathrm{d}s},\quad 0\le t\le T,\quad x\in\mathbb{R}. \tag{A18}$$

By construction, this is a (closed form) analytical solution as well as an entropy solution in the sense of Kružkov (Theorem 4.14 in Holden and Risebro (2015); see also Kružkov (1970)) to the hyperbolic balance law (A12). This fact can be checked by considering that the solution (A18) does not have discontinuities and satisfies the following weak formulation for any $\phi\in C_b^\infty$ with compact support:

$$\begin{aligned}&\int_{\mathbb{R}}\phi(x)\lambda_T(x)\mathrm{d}x-\int_{\mathbb{R}}\phi(x)\lambda_0(x)\mathrm{d}x\\&+\int_0^T\int_{\mathbb{R}}\frac{\mathrm{d}\phi(x)}{\mathrm{d}x}u_t\lambda_t(x)\mathrm{d}x\mathrm{d}t-\int_0^T\int_{\mathbb{R}}\phi(x)\left(R_t\lambda_t(x)+\frac{\sigma_t^2}{2}\left(\lambda_t(x)\right)^2\right)\mathrm{d}x\mathrm{d}t=0\end{aligned}. \tag{A19}$$

The hyperbolic balance law (A12) admits only one bounded and nonnegative entropy solution, which is (A18) because (A14) is an ODE whose uniqueness is guaranteed. Finally, $\mu$ and $\lambda$ are nonnegative.

□

***Proof of Proposition 2***

The proof is slightly long and is therefore divided into multiple steps.

**Step 1: Semidiscrete system: formulation**

We start from the following semidiscrete version of the superprocess model, whose limit is a desired weak solution. We also show that the Laplace functional of this semidiscrete model converges to the right-hand side of (14). The idea is qualitatively similar to that in Yoshioka (2026a), but ours is slightly different because it addresses a problem in an unbounded domain.

Set $h>0$ and $x_i=ih$ ($i\in\mathbb{Z}$). The semidiscrete model, which is considered a finite difference semidiscretization of (1), is presented as follows:

$$\mathrm{d}X_t^{(i)}=\left(u_t\frac{X_t^{(i+1)}-X_t^{(i)}}{h}-R_tX_t^{(i)}+f_t^{(i)}\right)\mathrm{d}t+\sigma_t\sqrt{X_t^{(i)}}\frac{\mathrm{d}B_t^{(i)}}{\sqrt{h}},\quad t>0,\quad i\in\mathbb{Z} \tag{A20}$$

subject to the initial condition $X_0^{(i)}\ge 0$ ($i\in\mathbb{Z}$), where $\left(B_t^{(i)}\right)_{t\ge 0}$ is a standard one-dimensional Brownian motion and all the elements of $\left\{B_t^{(i)}\right\}_{i\in\mathbb{Z}}$ are mutually independent and $f_t^{(i)}=f_t(x_i)$. Here, $X_t^{(i)}$ can be understood as a discretization of $h^{-1}\int_{(i-1)h}^{ih} X_t(\mathrm{d}x)$. We assume $h\sum_i \phi(x_i) X_0^{(i)} \underset{h\to 0}{\to} \int_{\mathbb{R}} \phi(x) X_0(\mathrm{d}x)$ for any $\phi \in C_b^\infty$.

The system (A20) consists of infinitely many SDEs but has a cascading nature such that $X_t^{(i)}$ can be found if $X_t^{(i+1)}$ is available. Therefore, we consider a truncated system such that $X_t^{(i)} \equiv 0$ ($|i| \geq N+1$) for each fixed $N \in \mathbb{N}$ and finally take $N \to +\infty$. In this case, the semidiscrete system becomes

$$\mathrm{d}X_t^{(i)} = \left( u_t \frac{X_t^{(i+1)}\mathbb{I}(i<N) - X_t^{(i)}}{h} - R_t X_t^{(i)} + f_t^{(i)} \right)\mathrm{d}t + \sigma_t \sqrt{X_t^{(i)}} \frac{\mathrm{d}B_t^{(i)}}{\sqrt{h}}, \quad t>0, \quad |i| \leq N. \tag{A21}$$

Here, $\mathbb{I}(i<N)=1$ if $i<N$ and $\mathbb{I}(i<N)=0$ if $i \geq N$, and the other indicator functions $\mathbb{I}(\cdot)$ are defined in the same way. This system of SDEs admits a unique pathwise solution that is nonnegative and continuous. Indeed, the boundedness of the coefficients $u_t, R_t, f_t^{(i)}, \sigma_t$ allows for the application of an argument analogous to that in pp.5-6 in Alfonsi (2015) along with Corollary 3.23 and Theorem 4.22 in Chapter 5 of Karatzas and Shreve (1991). Specifically, we have

$$(0 \leq)\ \mathbb{E}\left[ \sup_{t\leq T} \left( \sum_{|i|\leq N} X_t^{(i)} h \right)^2 \right], \quad \sup_{|i|\leq N, t\leq T} \mathbb{E}\left[ X_t^{(i)} \right] / (Nh) \leq K \tag{A22}$$

for some $K>0$ independent of $N$ (Proofs of Propositions 5 and 6 in Yoshioka (2026a)). Here, the first one comes from **Assumption 1** and (A34) shown in the next step; indeed, due to (A32) and Cauchy–Schwarz inequality, we have

$$\left( h \sum_{|i|\leq N} X_t^{(i)} \right)^2 \leq 3 \left( \left( \omega \sup_{x\in\mathbb{R}} f_t(x) \right)^2 + \left( h \sum_{|i|\leq N} X_0^{(i)} \right)^2 + \left( h \sum_{|i|\leq N} \int_0^t \sigma_s \sqrt{X_s^{(i)}} \frac{\mathrm{d}B_s^{(i)}}{\sqrt{h}} \right)^2 \right). \tag{A23}$$

We can evaluate the last term in (A23) by using Doob's inequality (Theorem 3.8(iv) in Karatzas and Shreve (1991)) combined with (A34):

$$\mathbb{E}\left[ \sup_{0\leq t\leq T} \left( h \sum_{|i|\leq N} \int_0^t \sigma_s \sqrt{X_s^{(i)}} \frac{\mathrm{d}B_s^{(i)}}{\sqrt{h}} \right)^2 \right] \leq 4\mathbb{E}\left[ \left( h \sum_{|i|\leq N} \int_0^T \sigma_s \sqrt{X_s^{(i)}} \frac{\mathrm{d}B_s^{(i)}}{\sqrt{h}} \right)^2 \right] \leq 4 \sup_{0\leq t\leq T} \left(\sigma_t^2\right) \mathbb{E}\left[ \int_{\mathbb{R}} Y_{N,T}(\mathrm{d}x) \right]. \tag{A24}$$

**Step 2: Semidiscrete system: Laplace functional**

The Laplace functional of the semidiscrete system (A21). For any nonnegative $\psi \in C_b^\infty$, we set $\psi_i = \psi(x_i)$. We estimate the Laplace functional of the semidiscrete system (A21) as follows:

$$\mathbb{E}\left[ \exp\left( -\sum_{|i|\leq N} \psi_i X_T^{(i)} h \right) \middle| \mathcal{F}_t \right] = \exp\left( -\mu_t^{(N,h)} - \sum_{|i|\leq N} \lambda_t^{(i,N,h)} X_t^{(i)} h \right), \quad 0 \leq t \leq T. \tag{A25}$$

We express the right-hand side of (A25) as $\exp\left(-U_t^{(N,h)}\right)$. Applying Itô's formula yields

$$
\begin{aligned}
&\mathrm{d}\left(\exp\left(-U_t^{(N,h)}\right)\right)\\
&=-\exp\left(-U_t^{(N,h)}\right)\left(\begin{array}{l}\dfrac{\mathrm{d}\mu_t^{(N,h)}}{\mathrm{d}t}+\sum_{|i|\le N}\lambda_t^{(i,N,h)}f_t h\\ +\sum_{|i|\le N}\left(\begin{array}{l}\dfrac{\mathrm{d}\lambda_t^{(i,N,h)}}{\mathrm{d}t}-u_t\dfrac{\lambda_t^{(i,N,h)}-\lambda_t^{(i-1,N,h)}\mathbb{I}(i>-N)}{h}\\ -R_t\lambda_t^{(i,N,h)}-\dfrac{\sigma_t^2}{2}\left(\lambda_t^{(i,N,h)}\right)^2\end{array}\right)X_t^{(i)}h\end{array}\right)\mathrm{d}t \qquad \text{(A26)}\\
&-\exp\left(-U_t^{(N,h)}\right)\sum_{|i|\le N}\sigma_t\sqrt{X_t^{(i)}}\sqrt{h}\mathrm{d}B_t^{(i)},
\end{aligned}
$$

where we used

$$
\sum_{|i|\le N}\lambda_t^{(i,N,h)}\frac{X_t^{(i+1)}\mathbb{I}(i<N)-X_t^{(i)}}{h}=-\sum_{|i|\le N}X_t^{(i)}\frac{\lambda_t^{(i,N,h)}-\lambda_t^{(i-1,N,h)}(i>-N)}{h}. \qquad \text{(A27)}
$$

Because the left-hand side of (A25) is a martingale, the drift term of (A26) must vanish, yielding the system of time-backward ODEs to find $\mu$ and $\lambda$:

$$
\frac{\mathrm{d}\lambda_t^{(i,N,h)}}{\mathrm{d}t}-u_t\frac{\lambda_t^{(i,N,h)}-\lambda_t^{(i-1,N,h)}(i>-N)}{h}-R_t\lambda_t^{(i,N,h)}-\frac{\sigma_t^2}{2}\left(\lambda_t^{(i,N,h)}\right)^2=0,\quad 0\le t<T,\ |i|\le N \qquad \text{(A28)}
$$

and

$$
\frac{\mathrm{d}\mu_t^{(N,h)}}{\mathrm{d}t}+\sum_{|i|\le N}\lambda_t^{(i,N,h)}f_t h=0,\quad 0\le t<T, \qquad \text{(A29)}
$$

where the terminal conditions are $\mu_T^{(N,h)}=0$ and $\lambda_T^{(i,N,h)}=\psi_i$.

As in Proof of Proposition 4 in Yoshioka (2026a), we obtain the following results for the ODE system (A28); the system (A28) admits a unique Carathéodory solution (XVIII. Existence and Uniqueness Theorem in Chapter 10 of Walter (1998)) such that $0\le\lambda_t^{(i,N,h)}\le\max_{x\in\mathbb{R}}\psi(x)$ ($0\le t\le T$, $|i|\le N$) and this bound is free from $h,N$; the summation in (A29) is bounded irrespective of $h,N$, and therefore, $\mu_t^{(N,h)}$ ($0\le t\le T$) is well defined; the total variation of $\lambda_t^{(i,N,h)}$, defined by $\mathrm{TV}^{(N,h)}(\lambda)=\sum_{i=-N+1}^{N}\left|\lambda_t^{(i,N,h)}-\lambda_t^{(i-1,N,h)}\right|$, is bounded irrespective of $t$ but depends on the product $L=Nh$, where $L$ is half of the domain length; and there exists a constant $C'>0$ depending on $L$ such that $h\sum_{i=-N}^{N}\left|\frac{\mathrm{d}\lambda_t^{(i,N,h)}}{\mathrm{d}t}\right|\le C'$, and hence $h\sum_{i=-N}^{N}\left|\lambda_t^{(i,N,h)}-\lambda_s^{(i,N,h)}\right|\le C'|t-s|$ ($0\le s,t\le T$).

Now, for each $k\in\mathbb{N}$, we choose $N=N_k=k\lfloor 1/h\rfloor$ ($\lfloor 1/h\rfloor$ the smallest natural number that is equal to or larger than $1/h$) and consider a localization argument. We consider a sequence $\{h_j\}_{j=1,2,3,\ldots}$ with $h_j\underset{j\to+\infty}{\to}0$. The total variation in the semidiscrete $\lambda_t^{(\cdot,N_k,h_j)}$ on $[-k,k]$ under this setting is denoted by $\mathrm{TV}_k$. First, for $k=1$, we have $\mathrm{TV}_1 < +\infty$, and by Helly's Theorem (Corollary A.10 in Holden and Risebro (2015)), there exists a subsequence $\left\{h_{1,j}\right\}_{j=1,2,3,...}$ of $\left\{h_j\right\}_{j=1,2,3,...}$ such that $\lambda_t^{\left(\cdot,N_k,h_{1,j}\right)}$ converges to a function integrable on $\left[-1,1\right]$ as $j \to +\infty$. Second, for $k=2$, we have $\mathrm{TV}_2 < +\infty$ for $\lambda_t^{\left(\cdot,N_k,h_{1,j}\right)}$. By Helly's Theorem, there exists a subsequence $\left\{h_{2,j}\right\}_{j=1,2,3,...}$ of $\left\{h_{1,j}\right\}_{j=1,2,3,...}$ such that it converges to a function integrable on $\left[-2,2\right]$. By repeating this procedure for $k=1,2,3...$, we can extract a subsequence $\left\{h_{j,j}\right\}_{j=1,2,3,...}$ of $\left\{h_j\right\}_{j=1,2,3,...}$ such that $\lambda_t^{\left(\cdot,N_k,h_{j,j}\right)}$ converges on arbitrary bounded compact intervals, implying that the limit of the sequence $\lambda_t^{\left(\cdot,N_k,h_{j,j}\right)}$ is locally integrable in $\mathbb{R}$.

Consequently, there exists a time-dependent function on $\mathbb{R}$, e.g., $\lambda_t\left(\cdot\right)$, with bounded variation in each compact interval in $\mathbb{R}$ such that

$$0=\int_{\mathbb{R}}\phi\left(x\right)\left(\frac{\partial\lambda_t\left(x\right)}{\partial t}-R_t\lambda_t\left(x\right)-\frac{1}{2}\left(\sigma_t\lambda_t\left(x\right)\right)^2\right)\mathrm{d}x+\int_{\mathbb{R}}u_t\frac{\mathrm{d}\phi\left(x\right)}{\mathrm{d}x}\lambda_t\left(x\right)\mathrm{d}x=0 \tag{A30}$$

for all $\phi \in C_b^\infty$ with compact support. This is identical to the weak formulation of (A12), which admits the unique entropy solution (A18). Therefore, the chosen subsequence $\lambda_t^{\left(\cdot,N_k,h_{j,j}\right)}$ converges to (A18). This convergence result can be extended to the original sequence $\lambda_t^{\left(\cdot,N_k,h_j\right)}$. Indeed, if the original sequence does not converge to $\lambda$ as a whole when $k,j \to +\infty$, then it should be possible to select a subsequence that converges to a limit $\bar{\lambda}$ that is different from $\lambda$. However, by the same procedure, $\bar{\lambda}$ would also be an entropy solution. Owing to the uniqueness of the entropy solutions, any chosen subsequence must necessarily converge to $\lambda$. Therefore, it is proven that not only a subsequence but also the entire sequence of difference solutions converges to the entropy solution $\lambda$.

**Step 3: Semidiscrete system: Tightness 1**

**Steps 3-4** follow the lines of Chapter 1 in Etheridge (2000) with adaptations to our problem. Fix $N \in \mathbb{N}$. We define the nonnegative measure

$$Y_{N,t}\left(\mathrm{d}x\right)=h\sum_{|i|\le N}X_t^{(i)}\delta\left(x-x_i\right) \tag{A31}$$

and show that it is tight in the sense of Definition 14 in Etheridge (2000). We have

$$\begin{aligned}\int_{\mathbb{R}}Y_{N,t}\left(\mathrm{d}x\right)&=h\sum_{|i|\le N}X_t^{(i)}\\&=h\sum_{|i|\le N}\int_0^t\left(u_s\frac{X_s^{(i+1)}\mathbb{I}\left(i<N\right)-X_s^{(i)}}{h}-R_sX_s^{(i)}+f_s^{(i)}\right)\mathrm{d}s+h\sum_{|i|\le N}\int_0^t\sigma_s\sqrt{X_s^{(i)}}\frac{\mathrm{d}B_s^{(i)}}{\sqrt{h}}+h\sum_{|i|\le N}X_0^{(i)}\\&=-u_sX_s^{(-N)}+h\sum_{|i|\le N}\int_0^t\left(-R_sX_s^{(i)}+f_s^{(i)}\right)\mathrm{d}s+h\sum_{|i|\le N}\int_0^t\sigma_s\sqrt{X_s^{(i)}}\frac{\mathrm{d}B_s^{(i)}}{\sqrt{h}}+h\sum_{|i|\le N}X_0^{(i)}\\&\le T\omega\sup_{x\in\mathbb{R}}f_t\left(x\right)+h\sum_{|i|\le N}\int_0^t\sigma_s\sqrt{X_s^{(i)}}\frac{\mathrm{d}B_s^{(i)}}{\sqrt{h}}+h\sum_{|i|\le N}X_0^{(i)},\end{aligned} \tag{A32}$$

and therefore

$$\int_{\mathbb{R}} Y_{N,t}(\mathrm{d}x) \le T\omega \sup_{0\le t\le T}\sup_{x\in\mathbb{R}} f_t(x) + h\sum_{|i|\le N}\int_0^t \sigma_s \sqrt{X_s^{(i)}}\frac{\mathrm{d}B_s^{(i)}}{\sqrt{h}} + h\sum_{|i|\le N} X_0^{(i)}, \tag{A33}$$

and hence

$$\mathbb{E}\left[\int_{\mathbb{R}} Y_{N,T}(\mathrm{d}x)\right] \le T\omega \sup_{0\le t\le T}\sup_{x\in\mathbb{R}} f_t(x) + h\sum_{|i|\le N} X_0^{(i)} \le C_T, \tag{A34}$$

where $C_T > 0$ is a constant independent from $N, h$ due to **Assumption 1**. Therefore, for any $K' > 0$, Markov's inequality (Theorem 5.11 in Klenke (2020)) shows

$$\Pr\left(\sup_{0\le t\le T}\int_{\mathbb{R}} Y_{N,t}(\mathrm{d}x) \ge K'\right) \le \frac{1}{K'}\mathbb{E}\left[\sup_{0\le t\le T}\int_{\mathbb{R}} Y_{N,t}(\mathrm{d}x)\right] \underset{K'\to+\infty}{\to} 0. \tag{A35}$$

Indeed, this follows from (A33) along with Doob's inequality (Theorem 3.8(iv) in Karatzas and Shreve (1991)):

$$\begin{aligned}\mathbb{E}\left[\sup_{0\le t\le T}\int_{\mathbb{R}} Y_{N,t}(\mathrm{d}x)\right] &\le T\omega \sup_{0\le t\le T}\sup_{x\in\mathbb{R}} f_t(x) + h\sum_{|i|\le N} X_0^{(i)} + \mathbb{E}\left[\sup_{0\le t\le T}\left|h\sum_{|i|\le N}\int_0^t \sigma_s\sqrt{X_s^{(i)}}\frac{\mathrm{d}B_s^{(i)}}{\sqrt{h}}\right|\right] \\ &\le T\omega \sup_{0\le t\le T}\sup_{x\in\mathbb{R}} f_t(x) + h\sum_{|i|\le N} X_0^{(i)} + 2\left(\mathbb{E}\left[h\sum_{|i|\le N}\int_0^T \sigma_s^2 X_s^{(i)}\mathrm{d}s\right]\right)^{\frac{1}{2}} \\ &< +\infty,\end{aligned} \tag{A36}$$

irrespective of $N$ and $h$.

Now, owing to (A35), by choosing $N = \lfloor 1/h^2 \rfloor$ and letting $h \to +0$, the right-hand side of (A35) is still $O(K'^{-1})$; hence, the process $\int_{\mathbb{R}} Y_{N,t}(\mathrm{d}x)$ is tight at each $t \in [0,T]$ in $\mathcal{M}_+(\overline{\mathbb{R}})$ ($\overline{\mathbb{R}}$ is a one-point compactification of $\mathbb{R}$, and $\mathcal{M}_+(\overline{\mathbb{R}})$ is the collection of all finite measures on $\overline{\mathbb{R}}$: item 1 of Proposition 1.19 in Etheridge (2000)). Additionally, the boundedness of the coefficients $u, f, R, \sigma$ indicates that the process $Y_{N,t}(\mathrm{d}x)$ satisfies the Aldous criterion (item 1 of Proposition 1.19 in Etheridge (2000)); thus, the tightness of $Y_{N,\cdot}(\mathrm{d}x)$ follows.

**Step 3: Semidiscrete system: Tightness 2**

Fix $N \in \mathbb{N}$. For any nonnegative $\phi \in C_b^\infty$, we obtain

$$\begin{aligned}&\int_{\mathbb{R}} \phi(x) Y_{N,t}(\mathrm{d}x) - \int_{\mathbb{R}} \phi(x) Y_{0,t}(\mathrm{d}x) \\ &= \int_0^t \left(\sum_{|i|\le N}\left(-u_t \frac{\phi_i - \phi_{i-1}\mathbb{I}(i > -N)}{h} - R_t\phi_i - \frac{\sigma_t^2}{2}\phi_i^2\right) X_t^{(i)} h + \sum_{|i|\le N} f_t^{(i)} h\right)\mathrm{d}s + \int_0^t \sum_{|i|\le N}\sigma_s\phi_i\sqrt{X_s^{(i)}}\sqrt{h}\mathrm{d}B_s^{(i)},\end{aligned} \tag{A37}$$

and the last term is a martingale with quadratic variation $\int_0^t \sigma_s^2 \sum_{|i|\le N}\phi_i^2 X_s^{(i)} h\mathrm{d}s$. By taking $N = \lfloor 1/h^2 \rfloor$ and letting $h \to +0$ in (A37), any limit point $X$ of $Y_N$ is a martingale given by (p.14 in Etheridge (2000)):

$$\int_{\mathbb{R}} \phi(x) X_t(\mathrm{d}x) - \int_{\mathbb{R}} \phi(x) X_0(\mathrm{d}x) - \int_0^t \int_{\mathbb{R}} \left\{ \left( -u_t \frac{\mathrm{d}\phi(x)}{\mathrm{d}x} - R_t \phi(x) \right) X_s(\mathrm{d}x) + f_t(x) \mathrm{d}x \right\} \mathrm{d}s \tag{A38}$$

with quadratic variation (13).

Finally, for any $\varepsilon > 0$, Markov's inequality shows

$$\Pr\left( \int_{|x| \geq L} X_t(\mathrm{d}x) > \varepsilon \right) \leq \frac{1}{\varepsilon} \hat{\mathbb{E}}\left[ \int_{|x| \geq L} X_t(\mathrm{d}x) \right] \underset{L \to +\infty}{\to} 0, \tag{A39}$$

indicating that the limit $X$ belongs to $\mathcal{M}_+$. Indeed, for $L > 1$, we have

$$\mathbb{E}\left[ \exp\left( -\int_{|x| \geq L} X_t(\mathrm{d}x) \right) \right] = \mathbb{E}\left[ \exp\left( -\int_{\mathbb{R}} \mathbb{I}(|x| \geq L) X_t(\mathrm{d}x) \right) \right] \geq \mathbb{E}\left[ \exp\left( -\int_{\mathbb{R}} \left(1 - \phi^{(L)}(x)\right) X_t(\mathrm{d}x) \right) \right], \tag{A40}$$

where $\phi^{(L)} \in C_b^{\infty}$ is such that $\phi^{(L)}(x) = 1$ for $|x| \leq L - L^{-1}$ and $\phi^{(L)}(x) = 0$ for $|x| \geq L$ and is increasing (resp., decreasing) for $-L \leq x \leq -L + L^{-1}$ (resp., $L - L^{-1} \leq x \leq L$ ). This $\phi^{(L)}$ is considered a smooth approximation of the indicator function $\mathbb{I}(|x| < L)$ from below.

By **Proposition 1**, with $\psi = 1 - \phi^{(L)}$, we obtain

$$\mathbb{E}\left[ \exp\left( \int_{\mathbb{R}} \left(1 - \phi^{(L)}(x)\right) X_t(\mathrm{d}x) \right) \right] = \exp\left( -\mu_0 - \int_{\mathbb{R}} \lambda_0(x) X_0(\mathrm{d}x) \right)\Bigg|_{\psi = 1 - \phi^{(L)}} \underset{L \to +\infty}{\to} 1 \tag{A41}$$

because of the formula (15) and **Assumption 1** ( $X_0$ and $f$ have compact support). This implies

$$\mathbb{E}\left[ \exp\left( -\int_{|x| \geq L} X_t(\mathrm{d}x) \right) \right] \underset{L \to +\infty}{\to} 1, \tag{A42}$$

and hence $\int_{|x| \geq L} X_t(\mathrm{d}x) \underset{L \to +\infty}{\to} 0$ in probability, leading to (A39).

**Step 4: Uniqueness**

We show that the weak solutions the superprocess model are unique in the sense of distribution, i.e., in the sense of the Laplace functional, as in Corollary 1.25 of Etheridge (2000). To show this, it suffices to recall **Proposition 1**, where we need to show that the following quantity is a martingale:

$$\begin{aligned} &\exp\left( -\int_{\mathbb{R}} \phi(x) X_t(\mathrm{d}x) \right) - \exp\left( -\int_{\mathbb{R}} \phi(x) X_0(\mathrm{d}x) \right) \\ &+ \int_0^t \exp\left( -\int_{\mathbb{R}} \phi(x) X_s(\mathrm{d}x) \right) \int_{\mathbb{R}} \left\{ \left( -u_t \frac{\mathrm{d}\phi(x)}{\mathrm{d}x} - R_t \phi(x) - \frac{\sigma_t^2}{2} (\phi(x))^2 \right) X_s(\mathrm{d}x) + f_t(x) \mathrm{d}x \right\} \mathrm{d}s, \end{aligned} \tag{A43}$$

but this follows from an argument analogous to **Step 3** by noting that

$$\begin{aligned} &\exp\left( -\int_{\mathbb{R}} \phi(x) Y_{N,t}(\mathrm{d}x) \right) - \exp\left( -\int_{\mathbb{R}} \phi(x) Y_{N,0}(\mathrm{d}x) \right) \\ &= -\int_0^t \exp\left( -\int_{\mathbb{R}} \phi(x) Y_{N,s}(\mathrm{d}x) \right) \left( \sum_{|i| \leq N} \left( -u_t \frac{\phi_i - \phi_{i-1} \mathbb{I}(i > -N)}{h} - R_t \phi_i - \frac{\sigma_t^2}{2} \phi_i^2 \right) X_t^{(i)} h + \sum_{|i| \leq N} f_t^{(i)} h \right) \mathrm{d}s \\ &- \int_0^t \exp\left( -\int_{\mathbb{R}} \phi(x) Y_{N,s}(\mathrm{d}x) \right) \sum_{|i| \leq N} \sigma_s \phi_i \sqrt{X_s^{(i)}} \sqrt{h} \mathrm{d}B_s^{(i)} \end{aligned} \tag{A44}$$

whose last term is a martingale with quadratic variation $\int_0^t \sigma_s^2 \exp\left( -2\int_{\mathbb{R}} \phi(x) Y_{N,s}(\mathrm{d}x) \right) \sum_{|i| \leq N} \phi_i^2 X_s^{(i)} h \mathrm{d}s$ .

Finally, by (14), iterating the conditional expectation backward in time determines the joint Laplace functional at any finitely many times; hence, the finite-dimensional distributions, and thus the law of $X$, are unique.

□

***Proof of Proposition 3***

The proof is due to a direct method leveraging the closed-form solutions.

**(a) Lower-bound case**

Fix a bounded and measurable mapping $\eta$ adapted to the filtration generated by $X$, and let $\mathbb{Q}=\mathbb{Q}_\eta$ be the induced probability measure. Because $\eta$ is bounded, applying (21) with $\phi\equiv 1$ and Gronwall's inequality yields $\mathbb{E}_{\mathbb{Q}}\left[\int_{\mathbb{R}} X_t(\mathrm{d}x)\right]\le C$ with a constant $C>0$ depending on $T, X_0, f$.

Let $\Phi$ be given by (39) with (34) and (35), where $G_t(x)$ is bounded, smooth in $x$, and absolutely continuous in $t$. This $\Phi$ solves the HJB equation (40). Applying (21) to the time-dependent test function $G_t(x)$, which is possible by the same argument at the beginning of **Proof of Proposition 1**, we obtain that the following quantity is a martingale under $\mathbb{Q}$:

$$\begin{aligned}&\Phi(t,X_t)-\Phi(0,X_0)\\&-\int_0^t\left(\frac{\mathrm{d}F_s}{\mathrm{d}s}+\int_{\mathbb{R}} f_s(x)G_s(x)\mathrm{d}x+\int_{\mathbb{R}}\left(\frac{\partial G_s(x)}{\partial s}-u_s\frac{\partial G_s(x)}{\partial x}-R_sG_s(x)+\sigma_s\eta_sG_s(x)\right)X_s(\mathrm{d}x)\right)\mathrm{d}s\end{aligned}. \tag{A45}$$

Its quadratic variation is given by $\int_0^t\sigma_s^2\int_{\mathbb{R}}(G_s(x))^2X_s(\mathrm{d}x)\mathrm{d}s$, which is integrable under $\mathbb{Q}$. Taking the expectation $\mathbb{E}_{\mathbb{Q}}$ of (A45) along with the terminal condition $\Phi(T,X_T)=\int_{\mathbb{R}}\psi(x)X_T(\mathrm{d}x)$ and the equations (34) and (35) yields

$$\mathbb{E}_{\mathbb{Q}}\left[\Phi(T,X_T)\right]=\Phi(0,X_0)+\mathbb{E}_{\mathbb{Q}}\left[\int_0^T\int_{\mathbb{R}}\left(\frac{\sigma_s^2}{2w}(G_s(x))^2+\sigma_s\eta_sG_s(x)\right)X_s(\mathrm{d}x)\mathrm{d}s\right]. \tag{A46}$$

Applying (22) to (A46) yields

$$\begin{aligned}&\mathbb{E}_{\mathbb{Q}}\left[\int_{\mathbb{R}}\psi(x)X_T(\mathrm{d}x)\right]+\frac{w}{2}\mathbb{E}_{\mathbb{Q}}\left[\int_0^T\int_{\mathbb{R}}(\eta_s(X_s,x))^2X_s(\mathrm{d}x)\mathrm{d}s\right]-\Phi(0,X_0)\\&=\frac{w}{2}\mathbb{E}_{\mathbb{Q}}\left[\int_0^T\int_{\mathbb{R}}\left(\frac{\sigma_s^2}{w^2}(G_s(x))^2+\frac{2\sigma_s\eta_s}{w}G_s(x)+(\eta_s(X_s,x))^2\right)X_s(\mathrm{d}x)\mathrm{d}s\right]\\&=\frac{w}{2}\mathbb{E}_{\mathbb{Q}}\left[\int_0^T\int_{\mathbb{R}}\left(\frac{\sigma_s}{w}G_s(x)+\eta_s(X_s,x)\right)^2X_s(\mathrm{d}x)\mathrm{d}s\right]\\&\ge 0.\end{aligned} \tag{A47}$$

The last inequality becomes equality if and only if the last expectation vanishes, yielding $\eta_t(X_t,x)=-\sigma_tG_t(x)/w$. Consequently, we obtain $\underline{\Phi}=\Phi(0,X_0)$ and this $\eta$ is optimal.

**(b) Upper-bound case**

Under (38), the denominator of $G_t(x)$ in (37) is positive, and this $G_t(x)$ is bounded. As for the case (a), we obtain

$$\begin{aligned}&\mathbb{E}_{\mathbb{Q}}\left[\int_{\mathbb{R}}\psi(x)X_T(\mathrm{d}x)\right]-\frac{w}{2}\mathbb{E}_{\mathbb{Q}}\left[\int_0^T\int_{\mathbb{R}}\left(\eta_s(X_s,x)\right)^2X_s(\mathrm{d}x)\mathrm{d}s\right]-\Phi(0,X_0)\\&=-\frac{w}{2}\mathbb{E}_{\mathbb{Q}}\left[\int_0^T\int_{\mathbb{R}}\left(\frac{\sigma_s}{w}G_s(x)-\eta_s(X_s,x)\right)^2X_s(\mathrm{d}x)\mathrm{d}s\right]\\&\leq 0.\end{aligned}\tag{A48}$$

The last inequality becomes equality if and only if the last expectation vanishes, yielding $\eta_t(X_t,x)=\sigma_tG_t(x)/w$. Consequently, we obtain $\bar{\Phi}=\Phi(0,X_0)$ and this $\eta$ is optimal.

□

***Proof of Proposition 4***

We only prove the (a) lower-bound case because the proof for the (b) upper-bound case is essentially the same except that we need to assume (38) and (47); the former ensures that the denominator of $G_t(x)$ is strictly positive, and the latter is used in (A51) below.

We set $\sup_{x\in\mathbb{R}}\psi(x)=\bar{\psi}$. Then, we have

$$0\leq G_t(x)\leq\psi\left(x-\int_t^Tu_s\mathrm{d}s\right)\leq\bar{\psi},\quad 0\leq t\leq T,\quad x\in\mathbb{R},\tag{A49}$$

yielding $|\eta_t^*|<+\infty$ ($0\leq t\leq T$). Now, we apply Dawson's Girsanov theorem (Theorem 7.9 in Etheridge (2000)) to show that $Z$ with $\eta=\eta^*$ is a martingale. To adapt Proof of Theorem 7.9 in Etheridge (2000) to our case, it suffices to show

$$\mathbb{E}_{\mathbb{Q}}\left[\sup_{0\leq t\leq T}\left(\int_{\mathbb{R}}X_t(\mathrm{d}x)\right)^2\right]<+\infty.\tag{A50}$$

To show this, under $\mathbb{Q}$, by (47), we have

$$0\leq\int_{\mathbb{R}}X_t(\mathrm{d}x)\leq\int_{\mathbb{R}}X_0(\mathrm{d}x)+\int_0^t\int_{\mathbb{R}}f_s(x)\mathrm{d}x\mathrm{d}s+M_{\mathbb{Q},t}(1),\tag{A51}$$

and hence

$$\sup_{0\leq t\leq T}\left(\int_{\mathbb{R}}X_t(\mathrm{d}x)\right)^2\leq 2\left\{\left(\int_{\mathbb{R}}X_0(\mathrm{d}x)\right)^2+\left(\int_0^T\int_{\mathbb{R}}f_s(x)\mathrm{d}x\mathrm{d}s\right)^2+\sup_{0\leq t\leq T}\left(M_{\mathbb{Q},t}(1)\right)^2\right\}.\tag{A52}$$

Doob's maximal inequality (Theorem 3.8(iv) in Karatzas and Shreve (1991)) shows

$$\mathbb{E}_{\mathbb{Q}}\left[\sup_{0\leq t\leq T}\left(M_{\mathbb{Q},t}(1)\right)^2\right]\leq 4\mathbb{E}_{\mathbb{Q}}\left[\left(M_{\mathbb{Q},T}(1)\right)^2\right]=4C'\int_0^t\hat{\mathbb{E}}_{\mathbb{Q}}\left[\int_{\mathbb{R}}X_s(\mathrm{d}x)\right]\mathrm{d}s\tag{A53}$$

with a constant $C'>0$. We also have

$$\mathbb{E}_{\mathbb{Q}}\left[\int_{\mathbb{R}}X_s(\mathrm{d}x)\right]\leq\int_{\mathbb{R}}X_0(\mathrm{d}x)+\int_0^t\int_{\mathbb{R}}f_s(x)\mathrm{d}x\mathrm{d}s<+\infty.\tag{A54}$$

Consequently, we arrive at the desired estimate (A50).

□

### A.2 Definitions of derivatives

Following Definitions 4.3 and 4.4 in Ocello (2025) (see also Definitions 2.1 and 2.3 in Martini (2023)), we present the definition of flat and intrinsic derivatives used in this study. First, for any continuous function $\phi:\mathcal{M}_+ \to \mathbb{R}$, a frat derivative, often called linear derivative as well, is defined as a mapping $\delta\phi:\mathbb{R}\times\mathcal{M}_+ \to \mathbb{R}$ such that

$$\phi(g)-\phi(g') = \int_0^1 \int_{\mathbb{R}} \delta\phi\left(x, mg+(1-m)g'\right)\left(g(\mathrm{d}x)-g'(\mathrm{d}x)\right)\mathrm{d}m \tag{A55}$$

if it is continuous for the product topology (continuous for any $x \in \mathbb{R}$ and weakly continuous for any $g \in \mathcal{M}_+$). Second, in the same setting, an intrinsic derivative of $D_x\phi:\mathbb{R}\times\mathcal{M}_+ \to \mathbb{R}$ at $x \in \mathbb{R}$ is defined as follows if $\delta\phi(x,g)$ is continuously differentiable at $x$ for any $g \in \mathcal{M}_+$:

$$D_x\phi(x,g) = \partial_x \delta\phi(x,g). \tag{A56}$$

Our $\Phi$ guessed in **Section 3** satisfies the regularity requirements for these derivatives for a.e. $0 \le t \le T$ due to the specific form (39).

### A.3 eDNA data

**Fig. 8** compares the empirical eDNA concentration and corresponding theoretical average concentration based on the fitted model for 2025. **Figs. 10 and 11** show sample paths and standard deviations of eDNA concentration with $\sigma = 2$ and $\sigma = \sqrt{10}$, respectively. eDNA concentrations and their variances are significantly larger than the case $\sigma = 1$.

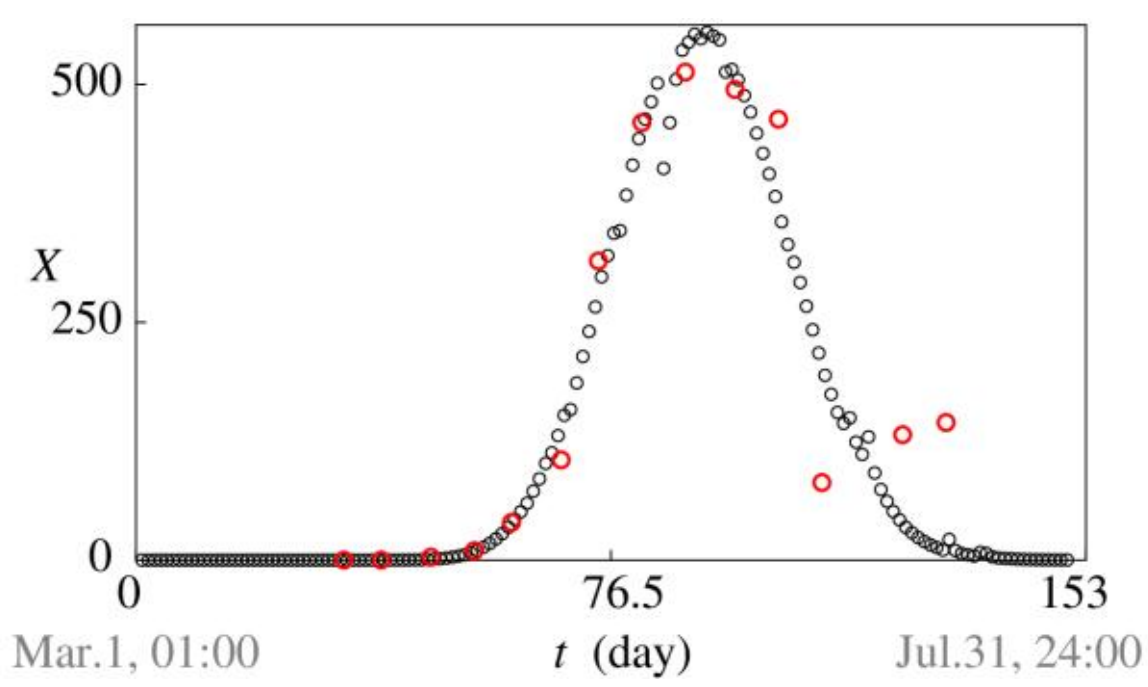

**Fig. 8** Comparison of eDNA concentrations $X$ (copies/mL) at the Kyu-Shita Bridge from the beginning of March to the end of July in 2025: empirical data (red) and theoretical average (black)

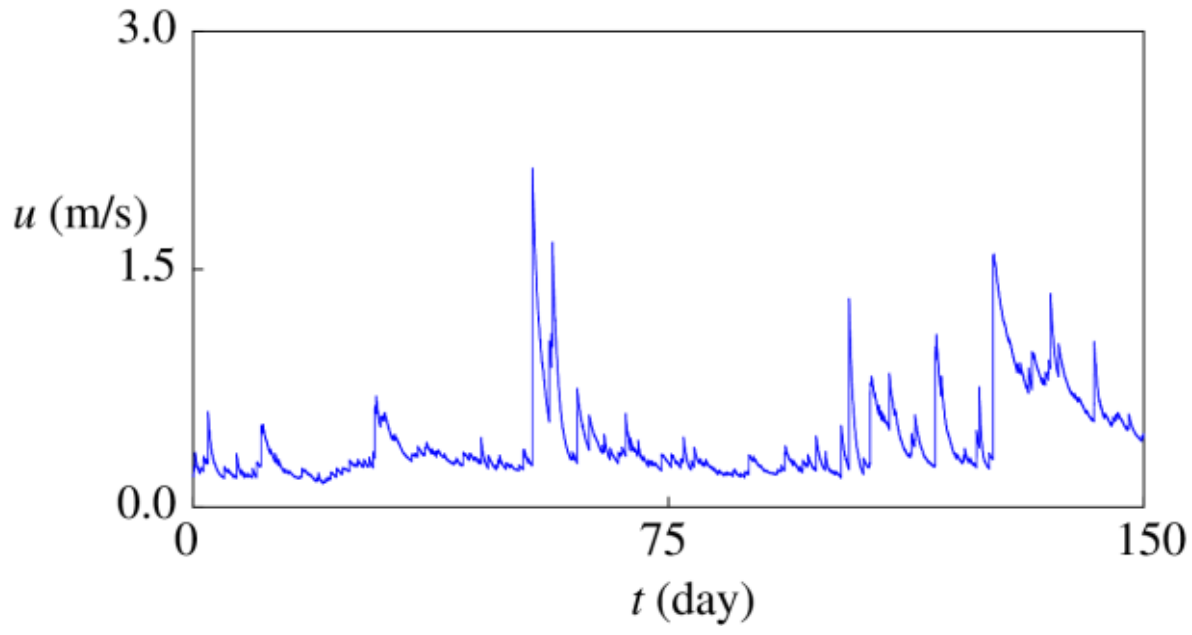

**Fig. 9** A sample path of the identified supOU process

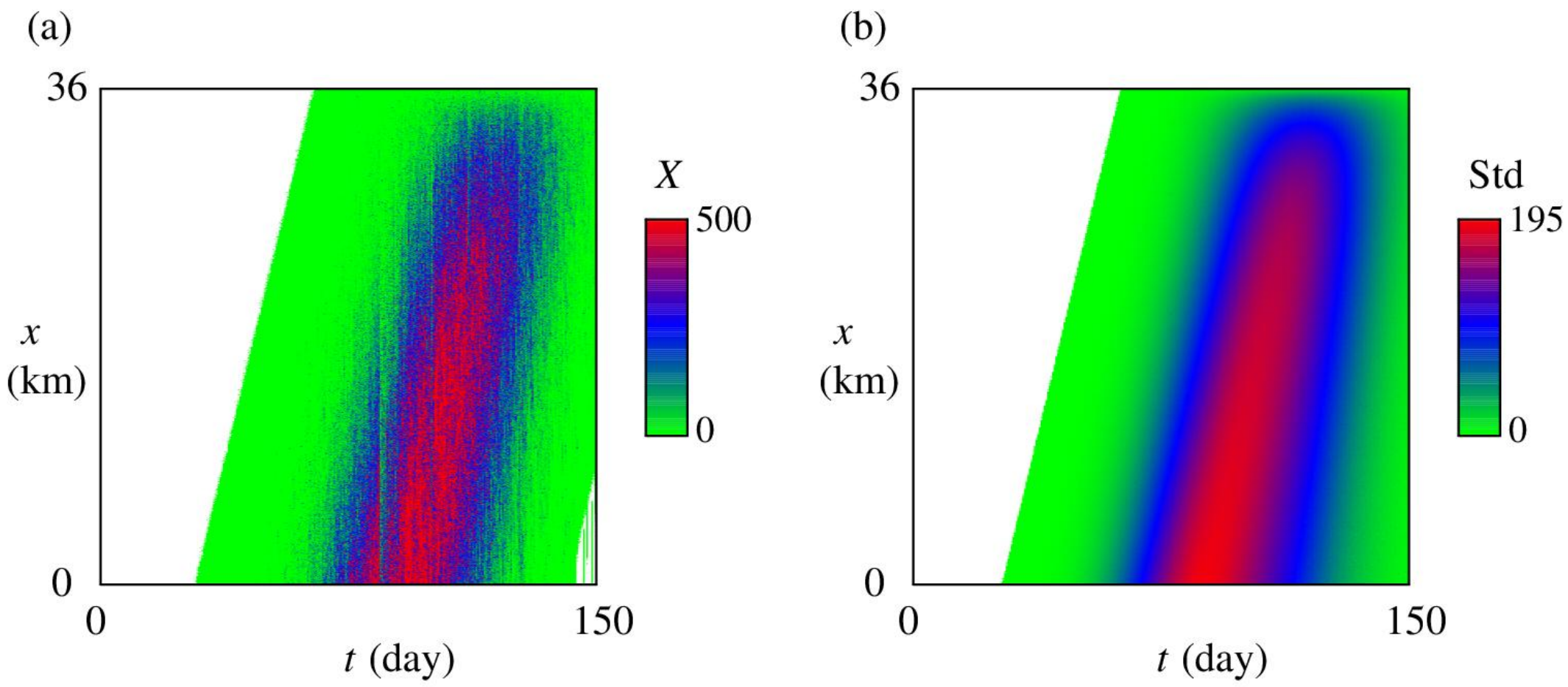

**Fig. 10** (a) Sample path of eDNA concentration $X$ (copies/mL) and (b) standard deviation (Std) of $X$ with $\sigma = 2$. The maximum was truncated by 500 for visualization ($\max X = 1958$ (copies/mL) in the space-time domain). The white area represents zero value

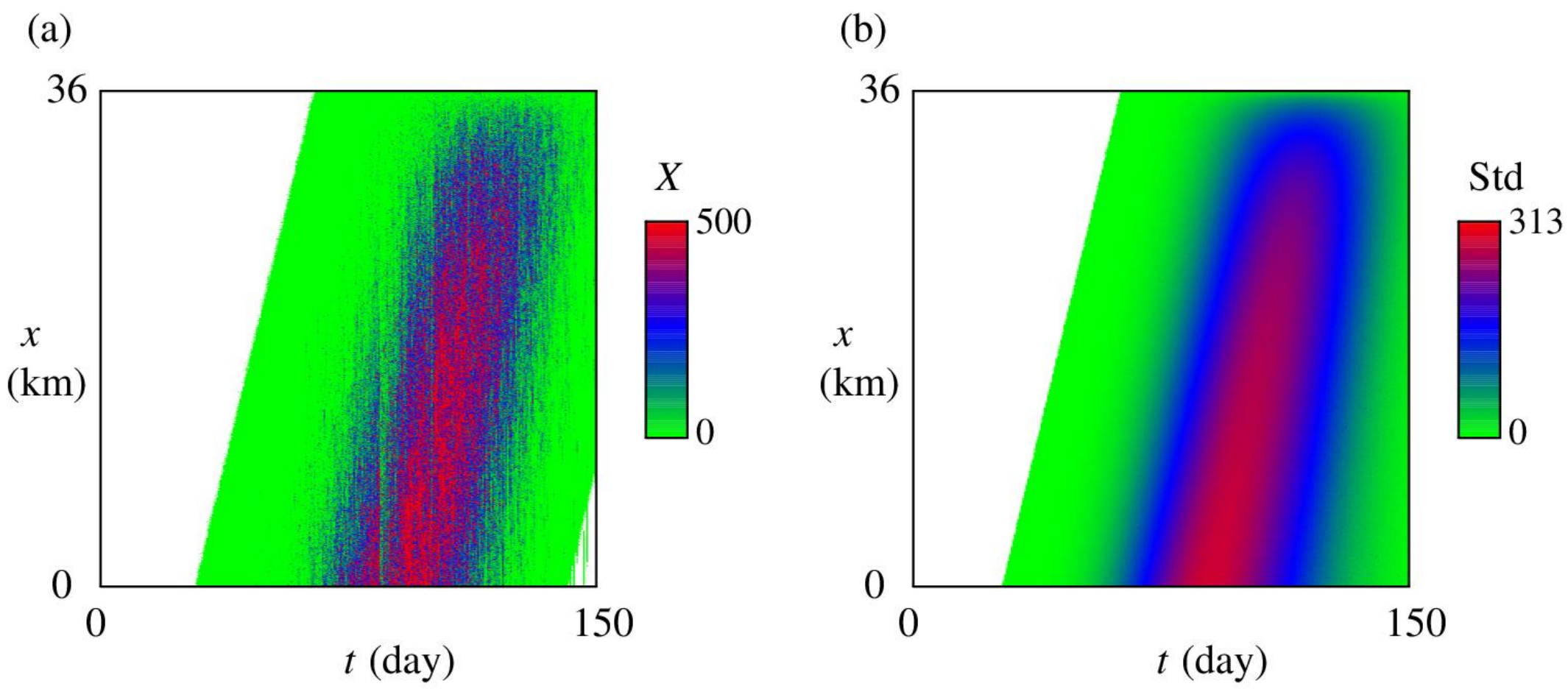

**Fig. 11** Same as **Fig. 10** but with $\sigma = \sqrt{10}$. The maximum was truncated by 500 for visualization ($\max X = 3416$ (copies/mL) in the space-time domain). The white area represents zero value

***Statements and Declarations***

**Funding:** This study was supported by the Japan Science and Technology Agency: PRESTO No. JPMJPR24KE and the Japan Society for the Promotion of Science: KAKENHI 26K00357.

**Competing Interests:** The authors have no relevant financial or nonfinancial interests to disclose.

**Data Availability:** The eDNA data are presented in **Appendix**. The other data will be made available upon reasonable request to the corresponding author.

**Declaration of Generative AI in Scientific Writing:** Claude (Anthropic) was used to assist in drafting parts of the text (the Introduction and some explanatory paragraphs); all content was verified and revised by the authors, who take full responsibility for the final manuscript.

**Author Contributions:**

**(Hidekazu Yoshioka)** Conceptualization, Methodology, Software, Formal analysis, Investigation, Resources, Data curation, Writing – original draft, Writing – review and editing, Visualization, Project administration, Funding acquisition

**(Motoh Tsujimura)** Formal analysis, Investigation, Writing – review and editing, Funding acquisition

**(Ayumi Hashiguchi)** Investigation, Writing – review and editing